\documentclass[11pt]{article}
\usepackage{amsmath,amssymb,amsthm,mathtools}
\usepackage[margin=2.5cm]{geometry}
\usepackage{multirow}
\usepackage{graphicx}
\usepackage[numbers,sort&compress]{natbib}
\newtheorem{theorem}{Theorem}[section]
\newtheorem{lemma}[theorem]{Lemma}
\newtheorem{proposition}[theorem]{Proposition}
\newtheorem{corollary}[theorem]{Corollary}
\theoremstyle{remark}

\newcommand{\E}{\mathbb{E}}
\newcommand{\Var}{\operatorname{Var}}
\newcommand{\Prob}{\mathbb{P}}
\newcommand{\R}{\mathbb{R}}
\newcommand{\dW}{d_{W}}

\newcommand{\dTV}{d_{TV}}
\newcommand{\bQ}{\bar Q}
\newcommand{\norm}[2]{\lVert #1\rVert_{#2}}

\title{Berry-Esséen bounds for Fractional Vasicek-Type Processes:\\
	Continuous and High-Frequency Regimes}
\author{Soukaina Douissi$^1$\\
	\small Cadi Ayyad University}
\date{\today}
\begin{document}
\maketitle

\footnotetext[1]{National School of Applied Sciences, Cadi Ayyad University, Marrakesh, Morocco. Email: s.douissi@uca.ac.ma}

\begin{abstract}
We study the estimation of the drift parameters $(a,b)$ of the fractional Vasicek process $dX_t=a(b-X_t)\,dt+dB^H_t$, driven by a fractional Brownian motion of Hurst index $H\in(0,1)$. The estimators are of moment type, built from the empirical mean and variance of the path, and are analysed in two regimes: continuous observation on $[0,T]$, and high-frequency observation at times $t_i=i\Delta_n$ with $\Delta_n\to0$ and $T_n=n\Delta_n\to\infty$. Using the Malliavin--Stein method on the second Wiener chaos, we prove consistency of both estimators for every $H$, and quantitative central limit theorems in the Wasserstein distance: for $H<\frac34$ for the estimator of $a$, with explicit rates depending on $H$, and for every $H\in(0,1)$ for the estimator of $b$, at the rate $T^{-H}$. In the discrete regime the bounds contain an additional discretization term of order $\Delta_n^{\min\{H+1/2,\,1\}}$, and for $H<\frac12$ the central limit theorem for the estimator of $b$ requires in addition $T_n^{1-2H}\Delta_n^{1+2H}\to0$. A short Monte Carlo study illustrates the results.

\medskip
\noindent\textbf{Keywords:} Fractional Vasicek process; parameter estimation; Berry--Ess\'een bounds; Wasserstein distance; Malliavin-Stein method; High-frequency observations.

\noindent\textbf{MSC 2020:} 60F05; 60G22; 62F12; 60H07.
\end{abstract}

\section{Introduction}\label{sec:intro}

Mean-reverting diffusions are the basic building block of many models in mathematical finance. The Vasicek model of the short interest rate \cite{Vasicek1977} is the Ornstein--Uhlenbeck equation $dX_t=a(b-X_t)\,dt+dW_t$, in which $a>0$ is the speed of mean reversion and $b$ the long-run level. Replacing the Brownian motion $W$ by a fractional Brownian motion (fBm) $B^H$ of Hurst index $H\in(0,1)$ \cite{MandelbrotVanNess1968} allows the noise to have memory: long memory when $H>\frac12$, as in the continuous-time long-memory volatility models of Comte and Renault~\cite{ComteRenault1998}, and rough, anti-persistent trajectories when $H<\frac12$, which is the regime reported for the log-volatility of financial assets by Gatheral et~al.~\cite{GatheralJaissonRosenbaum2018}. In this paper we consider the fractional Vasicek process
\begin{equation}\label{eq:intro-model}
X_0=0,\qquad dX_t=a\,(b-X_t)\,dt+dB^H_t,\qquad t\ge0,
\end{equation}
where $a>0$ and $b\in\mathbb R$ are unknown, $H\in(0,1)$ is known, and $X$ is observed either continuously on $[0,T]$ or at the discrete times $t_i=i\Delta_n$, $i=0,\dots,n$, with $\Delta_n\to0$ and $T_n:=n\Delta_n\to\infty$ (high-frequency regime). The process $X$ is asymptotically stationary and ergodic \cite{CheriditoKawaguchiMaejima2003}, and $b$ is its asymptotic mean.

Drift estimation for equations driven by fBm has been studied extensively. For the fractional Ornstein--Uhlenbeck (fOU) process, Hu and Nualart~\cite{HuNualart2010} studied a least-squares estimator for $H\ge\frac12$, proved its strong consistency and a central limit theorem for $H\in[\frac12,\frac34)$, and introduced an alternative, simulation-friendly estimator based on the ergodic theorem. Hu et~al.~\cite{HuNualartZhou2019} extended the study to the whole range $H\in(0,1)$, for continuous as well as discrete observations, and treated the case where the volatility is also unknown. Closer to the present setting, Mishura and Ralchenko~\cite{MishuraRalchenko2014} studied drift estimation in models driven by fBm from discrete observations, and Kukush et~al.~\cite{KukushMishuraRalchenko2017} proved the strong consistency of the drift estimators of the fOU process, for continuous and discrete observations and all $H\in(0,1)$. The estimators of the second type, which are also the ones used in this paper, are moment-type estimators based on the ergodic theorem, whose rates are governed by the covariance function of the stationary solution \cite{EsSebaiyViens2019}. For the Vasicek-type model~\eqref{eq:intro-model} with a continuous record, the least-squares estimators were analysed by Xiao and Yu~\cite{XiaoYu2019} for $H\ge\frac12$ and by Xiao and Yu~\cite{XiaoYu2019b} for $H<\frac12$, and the maximum likelihood estimator by Tanaka et~al.~\cite{TanakaXiaoYu2020}; Nourdin and Tran~\cite{NourdinTran2019} considered Hermite-driven Vasicek-type models and Alazemi et~al.~\cite{AlazemiDouissiEsSebaiy2019} studied discrete observations with a fixed step.

The Malliavin--Stein method of Nourdin and Peccati~\cite{NourdinPeccati2012}, together with the optimal fourth moment theorem \cite{NourdinPeccati2015}, has made it possible to obtain explicit rates of convergence in the central limit theorems for estimators of this type. Such Berry--Ess\'een bounds in the Wasserstein distance were obtained for second-moment estimators of stationary Gaussian processes and for the fOU process by Es-Sebaiy and Viens~\cite{EsSebaiyViens2019} and Douissi et~al.~\cite{DouissiEsSebaiyViens2019}, and in the least-squares setting by Chen et~al.~\cite{ChenKuangLi2020} and, for $H<\frac12$, by Chen and Li~\cite{ChenLi2021}. For observations at high frequency, Douissi et~al.~\cite{DouissiEsSebaiyKerchevNourdin2022} combined the Berry--Ess\'een bound of the continuous-time estimator with an estimate of the discretization error. This discretization error has recently been sharpened by Alazemi et~al.~\cite{AlazemiEsSebaiyAlForaih2026} and Es-Sebaiy and Chen~\cite{EsSebaiyChen2026}, who also treat the Kolmogorov and total variation distances. To our knowledge, however, the quantitative results available in the high-frequency regime concern estimators of a single drift parameter of a mean-zero model, whereas the Vasicek-type model~\eqref{eq:intro-model} requires the joint estimation of the speed $a$ and of the level $b$.

We construct moment-type estimators $(\hat a,\hat b)$ of $(a,b)$ from the empirical mean and the empirical variance of the path: $\hat b$ is the time average of $X$, which is a first-chaos functional of $B^H$, and $\hat a$ is obtained by inverting the relation $f_H(a)=\lim\bQ$ (almost surely), where $\bQ$ is the empirical variance and $f_H(a)=H\Gamma(2H)a^{-2H}$. This lives in the second Wiener chaos up to a negligible term, and our main results are the following.
\begin{itemize}
\item \emph{Continuous observation.} Both estimators are strongly consistent for every $H\in(0,1)$. For $H<\frac34$, we obtain explicit Berry--Ess\'een bounds in the Wasserstein distance for the empirical variance and for $\hat a$ (Theorem~\ref{thm:Q-cont}), and for every $H\in(0,1)$ we show that $T^{1-H}(\hat b_T-b)$ is at Wasserstein distance $O(T^{-H})$ of a centered Gaussian law (Theorem~\ref{thm:b-cont}).
\item \emph{High-frequency observation.} We prove consistency (Theorem~\ref{thm:consistency-discrete}) and Berry--Ess\'een bounds (Theorem~\ref{thm:BE-discrete}) in which the continuous-time rate is perturbed by a discretization term of order $\Delta_n^{\min\{H+1/2,\,1\}}$. For $H\ge\frac12$ no relation between $\Delta_n$ and $T_n$ is required for the central limit theorems; for $H<\frac12$ the central limit theorem for $\hat b_n$ needs $T_n^{1-2H}\Delta_n^{1+2H}\to0$ (Corollary~\ref{cor:CLT-discrete}). For $H\ge\frac12$ the discretization term is of the same order as the one obtained for the fOU process in Es-Sebaiy and Chen~\cite{EsSebaiyChen2026}.
\item \emph{Technical tools.} The estimator $\hat a$ requires uniform negative moments of the empirical variance. We prove them in the appendix with the Mandelbrot--Van Ness representation of the fBm and the Carbery--Wright inequality \cite{CarberyWright2001}, for the continuous and the discrete cases.
\end{itemize}
A short Monte Carlo study (Section~\ref{subsec:numerics}) illustrates the Gaussian approximation and the role of the slow rates near $H=\frac34$.

The paper is organized as follows. Section~\ref{sec:prelim} recalls the elements of Malliavin calculus on Wiener space that are used throughout. Section~\ref{sec:continuous} treats the continuous-time estimators, Section~\ref{sec:b-cont} the estimator of $b$, and Section~\ref{sec:discrete} the high-frequency observations, including the numerical illustration. The technical estimates on the covariance of the stationary solution and on the negative moments of the empirical variance are gathered in the Appendix.

\section{Preliminaries}\label{sec:prelim}

In this section we recall the elements of Malliavin calculus on Wiener space that are used throughout the paper; we refer to Nourdin and Peccati~\cite[Chapter~2]{NourdinPeccati2012} and Nualart~\cite[Chapter~1]{Nualart2006} for more details.

\subsection{Fractional Brownian motion and its Gaussian space}

Let $B^H=(B^H_t)_{t\in\mathbb R}$ be a two-sided fractional Brownian motion of Hurst index $H\in(0,1)$, that is, a centered Gaussian process with $B^H_0=0$ and covariance
\[
\E[B^H_tB^H_s]=R_H(t,s)=\tfrac12\bigl(|t|^{2H}+|s|^{2H}-|t-s|^{2H}\bigr).
\]
Let $\mathfrak H$ be the Hilbert space defined as the closure of the real-valued step functions on $\mathbb R$ with respect to the inner product $\langle\mathbf 1_{[0,t]},\mathbf 1_{[0,s]}\rangle_{\mathfrak H}=R_H(t,s)$. The map $\mathbf 1_{[0,t]}\mapsto B^H_t$ extends to a linear isometry $h\mapsto B^H(h)$ from $\mathfrak H$ onto the Gaussian space generated by $B^H$, so that $B^H=\{B^H(h):h\in\mathfrak H\}$ is an \emph{isonormal Gaussian process} over $\mathfrak H$:
\[
\E\bigl[B^H(h)B^H(g)\bigr]=\langle h,g\rangle_{\mathfrak H},\qquad h,g\in\mathfrak H .
\]
For deterministic kernels in $\mathfrak H$, $B^H(h)$ is the Wiener integral $\int_{\mathbb R}h\,dB^H$. It is convenient to also recall the Mandelbrot--Van Ness representation \cite{MandelbrotVanNess1968}
\begin{equation}\label{eq:MVN}
B^H_t=\int_{\mathbb R}K(t,u)\,dW_u,\quad t \in \R, \quad K(t,u)=c_H\bigl[(t-u)_+^{H-1/2}-(-u)_+^{H-1/2}\bigr],
\end{equation}
where $W$ is a two-sided standard Brownian motion and $c_H>0$ is such that $\Var B^H_1=1$.

\subsection{Wiener chaos and multiple integrals}

For $q\ge1$ let $\mathfrak H^{\otimes q}$ and $\mathfrak H^{\odot q}$ denote the $q$-th tensor product of $\mathfrak H$ and its symmetric subspace. The multiple Wiener--It\^o integral $I_q:\mathfrak H^{\odot q}\to L^2(\Omega)$ is defined by $I_q(h^{\otimes q})=H_q(B^H(h))$ for $\|h\|_{\mathfrak H}=1$, where $H_q$ is the $q$-th Hermite polynomial; $I_1(h)=B^H(h)$, and $I_q(f)=I_q(\tilde f)$ for $f\in\mathfrak H^{\otimes q}$, with $\tilde f$ the symmetrization of $f$. The $q$-th \emph{Wiener chaos} $\mathcal H_q$ is the closed linear span of $\{H_q(B^H(h)):\|h\|_{\mathfrak H}=1\}$, and $L^2(\Omega)=\bigoplus_{q\ge0}\mathcal H_q$ when the underlying $\sigma$-field is generated by $B^H$. We shall use the following properties.
\begin{itemize}
\item \textbf{Isometry.} For $f\in\mathfrak H^{\odot p}$ and $g\in\mathfrak H^{\odot q}$,
\begin{equation}\label{eq:isometry}
\E\bigl[I_p(f)I_q(g)\bigr]=\begin{cases}p!\,\langle f,g\rangle_{\mathfrak H^{\otimes p}}&\text{if }p=q,\\ 0&\text{otherwise.}\end{cases}
\end{equation}
In particular, elements of $\mathcal H_q$ with $q\ge1$ are centered, and chaoses of different orders are orthogonal.
\item \textbf{Contractions.} For $f\in\mathfrak H^{\odot p}$, $g\in\mathfrak H^{\odot q}$ and $0\le r\le p\wedge q$, the $r$-th contraction $f\otimes_r g\in\mathfrak H^{\otimes(p+q-2r)}$ is the function of $p+q-2r$ variables obtained by pairing the last $r$ variables of $f$ with the last $r$ variables of $g$ through the inner product of $\mathfrak H$:
\[
(f\otimes_r g)(t_1,\dots,t_{p-r},s_1,\dots,s_{q-r})=\Bigl\langle f(t_1,\dots,t_{p-r},\cdot),\,g(s_1,\dots,s_{q-r},\cdot)\Bigr\rangle_{\mathfrak H^{\otimes r}} ,
\]
with the convention $f\otimes_0 g=f\otimes g$. For $H>\frac12$, where $\langle\mathbf 1_{[0,t]},\mathbf 1_{[0,s]}\rangle_{\mathfrak H}=\alpha_H\int_0^t\!\!\int_0^s|u-v|^{2H-2}du\,dv$ with $\alpha_H=H(2H-1)$, this reads
\[
(f\otimes_r g)(t,s)=\alpha_H^{\,r}\int_{\mathbb R^{2r}}f(t,u)\,g(s,v)\prod_{\ell=1}^r|u_\ell-v_\ell|^{2H-2}\,du\,dv ,
\]
where $t=(t_1,\dots,t_{p-r})$, $s=(s_1,\dots,s_{q-r})$, $u=(u_1,\dots,u_r)$ and $v=(v_1,\dots,v_r)$. For $H\le\frac12$ the same notation is understood through the inner product of $\mathfrak H$.
We write $f\,\widetilde\otimes_r g$ for the symmetrization of $f\otimes_r g$.
\item \textbf{Product formula.} For $f\in\mathfrak H^{\odot p}$ and $g\in\mathfrak H^{\odot q}$,
\begin{equation}\label{eq:product}
I_p(f)\,I_q(g)=\sum_{r=0}^{p\wedge q}r!\binom pr\binom qr\,I_{p+q-2r}\bigl(f\,\widetilde\otimes_r g\bigr).
\end{equation}
The case used repeatedly below is $p=q=1$: for $f,g\in\mathfrak H$,
\begin{equation}\label{eq:product11}
I_1(f)\,I_1(g)=I_2\bigl(f\,\widetilde\otimes\, g\bigr)+\langle f,g\rangle_{\mathfrak H},\qquad\text{and in particular}\qquad I_1(f)^2=I_2(f\otimes f)+\|f\|_{\mathfrak H}^2 .
\end{equation}
Consequently, products of two Gaussian variables belong to $\mathcal H_0\oplus\mathcal H_2$, and a quadratic functional of a Gaussian process is a second-chaos element up to its mean.
\item \textbf{Hypercontractivity.} For every $q\ge1$ and every $F\in\bigoplus_{l=0}^q\mathcal H_l$, all the $L^p$ norms are equivalent to the $L^2$ norm: for every $p\ge2$ there exists a constant $c_{p,q}$, depending only on $p$ and $q$, such that
\begin{equation}\label{eq:hypercontractivity}
\bigl(\E|F|^p\bigr)^{1/p}\le c_{p,q}\,\bigl(\E|F|^2\bigr)^{1/2}.
\end{equation}

\end{itemize}

\subsection{Distances and normal approximation}

For random variables $U,V$ with finite first moment, the Wasserstein distance is
\[
d_{\mathrm W}(U,V)=\sup_{\phi\in\mathrm{Lip}(1)}\bigl|\E\phi(U)-\E\phi(V)\bigr|,
\]
where $\mathrm{Lip}(1)$ is the set of $1$-Lipschitz functions; the total variation and Kolmogorov distances are denoted by $d_{\mathrm{TV}}$ and $d_{\mathrm{Kol}}$, and $d_{\mathrm{Kol}}\le d_{\mathrm{TV}}$. Two elementary facts are used throughout. First, for any $U,V$ on the same probability space and any $N$,
\begin{equation}\label{eq:coupling-prelim}
d_{\mathrm W}(U+V,N)\le d_{\mathrm W}(U,N)+\E|V| ,
\end{equation}

\subsection{Cumulants and the optimal fourth moment theorem}

The third and fourth cumulants of a random variable $F$ with finite fourth moment are :
\[
\kappa_3(F)=\E F^3-3\E F^2\,\E F+2(\E F)^3,\qquad
\kappa_4(F)=\E F^4-4\E F\,\E F^3-3(\E F^2)^2+12(\E F)^2\E F^2-6(\E F)^4 ;
\]
for $\E F=0$, they reduce to $\kappa_3(F)=\E F^3$ and $\kappa_4(F)=\E F^4-3(\E F^2)^2$. \\
For a second-chaos random variable $F=I_2(h)$, with $h\in\mathfrak H^{\odot2}$ \cite[see][]{NourdinPeccati2010, BiermeBonamiNourdinPeccati2012},
\begin{equation}\label{eq:cumulants}
\kappa_3(F)=8\,\langle h,h\otimes_1h\rangle_{\mathfrak H^{\otimes2}},\qquad
|\kappa_4(F)|\le48\,\|h\otimes_1h\|^2_{\mathfrak H^{\otimes2}} .
\end{equation}
The optimal fourth moment theorem \cite{NourdinPeccati2015} states that, if $F_n$ lies in a fixed Wiener chaos $\mathcal H_q$, $\E F_n=0$, $\E F_n^2=1$ and $F_n$ converges in law to $\mathcal N(0,1)$, then there exist constants $0<c<C$, depending on the sequence but not on $n$, such that
\begin{equation}\label{eq:ofmt}
c\,\max\{|\kappa_3(F_n)|,\kappa_4(F_n)\}\le d_{\mathrm{TV}}\bigl(F_n,\mathcal N(0,1)\bigr)\le C\,\max\{|\kappa_3(F_n)|,\kappa_4(F_n)\}.
\end{equation}
By Douissi et~al.~\cite[Remark~2.2]{DouissiEsSebaiyKerchevNourdin2022} the same two-sided estimate holds with $d_{\mathrm W}$ in place of $d_{\mathrm{TV}}$.

\subsection{The Carbery--Wright inequality}

Let $Q$ be a polynomial of degree at most $d$ in $n$ independent standard Gaussian variables $X=(X_1,\dots,X_n)$. The Carbery--Wright inequality \cite{CarberyWright2001} states that there exists a universal constant $c>0$ (independent of $n$ and of $Q$) such that, for every $\varepsilon>0$,
\begin{equation}\label{eq:CW}
\Prob\bigl(|Q(X)|\le\varepsilon\bigr)\le c\,d\,\frac{\varepsilon^{1/d}}{\bigl(\E[Q(X)^2]\bigr)^{1/(2d)}} .
\end{equation}


\section{Continuous time observations}\label{sec:continuous}

In this paper, we deal with the drift parameter estimation problem for the fractional Vasicek process $X$, defined as the unique solution to the following stochastic differential equation
\begin{equation}
dX_t = a(b - X_t)dt + dB^H_t, \quad t \geqslant 0, \quad X_0 = 0, \tag{1}
\end{equation}
where $B^H$ is a fractional Brownian motion of Hurst parameter $H \in (0, 1)$. The Hurst parameter is supposed to be known whereas the parameters in (1) $a > 0$ and $b \in \mathbb{R}$ are unknown. We aim to estimate those parameters based on continuous and discrete time observation of the process $X$.

The process $X$ of equation (1) has the following explicit form
\begin{equation}
X_t = b(1 - e^{-at}) + \int_0^t e^{-a(t-u)} dB^H_u, \tag{2}
\end{equation}
where the integral with respect to $B^H$ is understood in the Wiener sense. When a complete path of the process $X$ over a large finite interval $[0, T]$ is observable, $T > 0$, we use the following estimators for $(a, b)$
\begin{equation}
\hat a_T := \left(\frac{\bar Q(T)}{H\Gamma(2H)}\right)^{-1/(2H)} \quad \text{and} \quad \hat b_T := \frac{1}{T}\int_0^T X_t dt. \tag{3}
\end{equation}
where
\begin{equation}
\bar Q(T) := \frac{1}{T}\int_0^T X_t^2 dt - \left(\frac{1}{T}\int_0^T X_t dt\right)^2. \tag{4}
\end{equation}
In the next step we will prove the strong consistency of $\hat a_T$ and $\hat b_T$ given in (3) for all $H \in (0, 1)$.

\subsection{Strong consistency}

\begin{theorem}\label{thm:consistency}
	Let $X$ be given by (1) and (2), where $B^H$ is a fractional Brownian motion of Hurst parameter $H \in (0, 1)$, and consider the estimators given in (3), then we have
	\begin{equation}
	(\hat a_T, \hat b_T) \to (a, b) \tag{5}
	\end{equation}
	almost surely as $T \to +\infty$.
\end{theorem}

\begin{proof}
	We can write the process $X$ in (1) and (2) as follows
	\begin{equation}
	X_t = b + Z_t^H + Y_t^H, \tag{6}
	\end{equation}
	where
	\begin{equation}
	Z_t^H := e^{-at}\int_{-\infty}^t e^{as} dB^H_s \quad \text{and} \quad Y_t^H := -e^{-at}(b + Z_0^H), \quad t \geqslant 0. \tag{7}
	\end{equation}
	Therefore, since the process $Z^H$ is an ergodic stationary Gaussian process with covariance $\rho_H(t) := \mathbb{E}[Z_t^H Z_0^H]$, by the Ergodic theorem, we have $\hat b_T \to b$, a.s.\ as $T \to +\infty$. Moreover, by the Ergodic theorem, the following convergence holds
	\begin{align*}
	\frac{1}{T}\int_0^T X_t^2 dt &= \frac{1}{T}\int_0^T (b + Z_t^H + Y_t^H)^2 dt \\
	&= b^2 + \frac{2b}{T}\int_0^T Z_t^H dt + \frac{1}{T}\int_0^T (Z_t^H)^2 dt + \frac{2b}{T}\int_0^T Y_t^H dt \\
	&\quad + \frac{2}{T}\int_0^T Z_t^H Y_t^H dt + \frac{1}{T}\int_0^T (Y_t^H)^2 dt \\
	&\longrightarrow b^2 + \mathbb{E}[(Z_0^H)^2] = b^2 + f_H(a). \tag{8}
	\end{align*}
	a.s.\ as $T \to +\infty$, where $f_H(a) := \frac{H\Gamma(2H)}{a^{2H}}$, therefore the strong consistency of $\bar Q_T$ and $\hat a_T$ defined in (4) and (3) respectively follows.
\end{proof}

\subsection{Asymptotics of the estimator $\hat a_T$}

The estimator $\hat a_T$ can be written as follows, $\hat a_T = f_H^{-1}(\bar Q_T)$, where $f_H(x) = \frac{H\Gamma(2H)}{x^{2H}}$. On the other hand, using the decomposition (6) of the process $X$, we get
\begin{equation}
\bar Q_T = \frac{1}{T}\int_0^T (Z_t^H)^2 dt + R_H(T), \tag{9}
\end{equation}
where
\begin{equation}
\begin{split}
R_H(T) := {}& \frac{2}{T}\int_0^T Z_t^H Y_t^H dt + \frac{1}{T}\int_0^T (Y_t^H)^2 dt - \left(\frac{1}{T}\int_0^T Z_t^H dt\right)^2 \\
& - \frac{2}{T^2}\int_0^T\int_0^T Z_t^H Y_s^H \,dt\,ds - \frac{1}{T^2}\left(\int_0^T Y_t^H dt\right)^2.
\end{split}
\tag{10}
\end{equation}

\begin{proposition}
	Consider $R_H(T)$ defined in (10), then we have
	\[
	\mathbb{E}[|R_H(T)|] \leqslant C \times
	\begin{cases}
	T^{-1} & \text{if } 0 < H \leqslant \frac{1}{2}, \\
	T^{2H-2} & \text{if } \frac{1}{2} \leqslant H < 1.
	\end{cases}
	\]
\end{proposition}

\begin{proof}
	\begin{align*}
	\mathbb{E}[R_H(T)] &= \frac{2}{T}\int_0^T \mathbb{E}[Z_t^H Y_t^H] dt + \frac{1}{T}\int_0^T \mathbb{E}[(Y_t^H)^2] dt - \frac{1}{T^2}\int_0^T\int_0^T \mathbb{E}[Z_t^H Z_s^H] dt\,ds \\
	&\quad - \frac{2}{T^2}\int_0^T\int_0^T \mathbb{E}[Z_t^H Y_s^H] dt\,ds - \frac{1}{T^2}\int_0^T\int_0^T \mathbb{E}[Y_t^H Y_s^H] dt\,ds.
	\end{align*}
	Moreover
	\begin{itemize}
		\item $\dfrac{2}{T}\displaystyle\int_0^T \mathbb{E}[Z_t^H Y_t^H] dt = -\dfrac{2}{T}\displaystyle\int_0^T \rho_H(t) e^{-at} dt \lesssim \dfrac{1}{T}$.
		\item $\dfrac{1}{T}\displaystyle\int_0^T \mathbb{E}[(Y_t^H)^2] dt = \dfrac{b^2 + f_H(a)}{T}\displaystyle\int_0^T e^{-2at} dt \lesssim \dfrac{1}{T}$.
		\item
		\begin{align*}
		\frac{1}{T^2}\int_0^T\int_0^T \mathbb{E}[Z_t^H Z_s^H] dt\,ds
		&= \frac{1}{T^2}\int_0^T\int_0^T \rho_H(t-s) dt\,ds
		= \frac{2}{T^2}\int_0^T\int_0^t \rho_H(u)\,du\,dt \\
		&= \frac{2}{T^2}\int_0^T\int_u^T \rho_H(u)\,du\,dt
		= \frac{2}{T}\int_0^T \left(1 - \frac{u}{T}\right)\rho_H(u)\,du \\
		&\lesssim
		\begin{cases}
		T^{-1} & \text{if } 0 < H < 1/2, \\
		T^{2H-2} & \text{if } 1/2 \leqslant H < 1.
		\end{cases}
		\end{align*}
		\item
		\begin{align*}
		\frac{2}{T^2}\int_0^T\int_0^T \mathbb{E}[Z_t^H Y_s^H] dt\,ds
		&= \frac{2}{T^2}\int_0^T e^{-as} ds \int_0^T \rho_H(t) dt
		\leqslant \frac{2}{aT^2}\int_0^T \rho_H(t) dt \\
		&\lesssim
		\begin{cases}
		T^{-2} & \text{if } 0 < H < 1/2, \\
		T^{2H-3} & \text{if } 1/2 \leqslant H < 1.
		\end{cases}
		\end{align*}
		\item $\dfrac{1}{T^2}\displaystyle\int_0^T\int_0^T \mathbb{E}[Y_t^H Y_s^H] dt\,ds = \dfrac{b^2+f_H(a)}{T^2}\left(\displaystyle\int_0^T e^{-at}dt\right)^2 \lesssim \dfrac{1}{T^2}$.
	\end{itemize}
	Consequently, using the decomposition (9) of $\bar Q_T$ and the fact that $Z^H$ is a stationary process along with Proposition 2, the following Lemma holds.
\end{proof}

\begin{lemma}
	The estimator $\bar Q_T$ is asymptotically unbiased for $f_H(a) = \frac{H\Gamma(2H)}{a^{2H}}$ for all $H \in (0, 1)$. More precisely, we have the following estimate for $T$ large
	\[
	|\mathbb{E}[\bar Q_T] - f_H(a)| \leqslant C \times
	\begin{cases}
	T^{-1} & \text{if } 0 < H \leqslant \frac{1}{2}, \\
	T^{2H-2} & \text{if } \frac{1}{2} \leqslant H < 1.
	\end{cases}
	\]
\end{lemma}

We have $Z_t^H = I_1(\varepsilon_t)$, where $\varepsilon_t(\cdot) := e^{-a(t-\cdot)} \mathbf{1}_{(-\infty,t]}$, therefore using the fact that $Z^H$ is stationary along with the product formula, we get
\[
\bar Q_T = I_2(f_T) + f_H(a) + R_H(T)
\]
with $f_T = \frac{1}{T}\int_0^T \varepsilon_t^{\otimes 2} dt$. It follows that
\[
\sqrt{T}\left(\bar Q_T - f_H(a)\right) = I_2(\sqrt{T} f_T) + \sqrt{T} R_H(T).
\]
Let us denote $v_T := \mathbb{E}[F_T^2]$, then
\begin{align*}
v_T &= 2\|f_T\|^2_{L^2(\mathbb{R}^2)} = \frac{2}{T}\int_0^T\int_0^T \langle \varepsilon_t^{\otimes 2}, \varepsilon_s^{\otimes 2}\rangle\,dt\,ds = \frac{2}{T}\int_0^T\int_0^T (\langle \varepsilon_t, \varepsilon_s\rangle)^2\,dt\,ds \\
&= \frac{2}{T}\int_0^T\int_0^T \left(\mathbb{E}[Z_t^H Z_s^H]\right)^2 dt\,ds = \frac{2}{T}\int_0^T\int_0^T \rho_H^2(t-s)\,dt\,ds = \frac{4}{T}\int_0^T\int_u^T \rho_H^2(u)\,dt\,du \\
&= 4\int_0^{+\infty}\left(1 - \frac{u}{T}\right)_+\rho_H^2(u)\,du
= 4\int_0^1\left(1-\frac{u}{T}\right)\rho_H^2(u)\,du + 4\int_1^T\left(1-\frac{u}{T}\right)\rho_H^2(u)\,du.
\end{align*}
Then, by Lebesgue dominated convergence theorem and using the fact that $|\rho_H(t)| = O(t^{2H-2})$ for large $t$, the following proposition follows.

\begin{proposition}
	$v_T \to \sigma_H^2 = 4\int_0^{+\infty}\rho_H^2(u)\,du$, as $T \to +\infty$.
\end{proposition}

Note that, throughout the paper, the symbol $\lesssim$ means that we omit multiplicative universal constants. The following theorem investigates the rate of convergence in law of $F_T$ towards a normal random variable when $T$ goes to infinity.

\begin{proposition}
	Assume that $0 < H \leqslant \frac{3}{4}$ and $N \sim \mathcal{N}(0,1)$, then
	\[
	\dTV(F_T, N) \lesssim
	\begin{cases}
	T^{-1/2} & \text{if } H \in \left(0, \frac{2}{3}\right) \\
	T^{-1/2}\log(T)^2 & \text{if } H = \frac{2}{3} \\
	T^{6H-9/2} & \text{if } H \in \left(\frac{2}{3}, \frac{3}{4}\right) \\
	\log(T)^{-3/2} & \text{if } H = \frac{3}{4}.
	\end{cases}
	\]
\end{proposition}

\begin{proof}
	By the optimal fourth moment theorem \cite{NourdinPeccati2015}, see \eqref{eq:ofmt}, and since $F_T$ belongs to the second Wiener chaos, we have
	\[
	\dTV\left(\frac{F_T}{\sqrt{v_T}}, \mathcal{N}(0,1)\right) \lesssim \left(\frac{|\kappa_3(F_T)|}{\mathbb{E}[F_T^2]^{3/2}}, \frac{|\kappa_4(F_T)|}{\mathbb{E}[F_T^2]^{2}}\right).
	\]
	On the other hand, see Appendix, one can prove that
	\[
	|\kappa_3(F_T)| \lesssim \frac{8}{\sqrt{T}}\left(\int_{-T}^T |\rho_H(t)|^{3/2} dt\right)^2 \quad \text{and} \quad |\kappa_4(F_T)| \lesssim \frac{48}{T}\left(\int_{-T}^T |\rho_H(t)|^{4/3} dt\right)^3.
	\]
\end{proof}

\begin{lemma}
	Assume that $0 < H < \frac{3}{4}$, then
	\begin{equation}
	|v_T - \sigma_H^2| \leqslant C
	\begin{cases}
	T^{-1} & \text{if } 0 < H \leqslant \frac{1}{2}, \\
	T^{4H-3} & \text{if } \frac{1}{2} < H < \frac{3}{4}.
	\end{cases} \tag{11}
	\end{equation}
\end{lemma}

\begin{proof}
	\begin{align*}
	|v_T - \sigma_H^2| &= \left|v_T - 4\int_0^{+\infty}\rho_H^2(u)\,du\right| \\
	&= \left|4\int_0^{+\infty}\left(1-\frac{u}{T}\right)_+ \rho_H^2(u)\,du - 4\int_0^{+\infty}\rho_H^2(u)\,du\right| \\
	&= \left|4\int_0^T\left(1-\frac{u}{T}\right)\rho_H^2(u)\,du - 4\int_0^{+\infty}\rho_H^2(u)\,du\right| \\
	&= \left|4\int_T^{+\infty}\rho_H^2(u)\,du - 4\int_0^T \frac{u}{T}\rho_H^2(u)\,du\right| \\
	&\leqslant C
	\begin{cases}
	T^{-1} & \text{if } 0 < H \leqslant \frac{1}{2}, \\
	T^{4H-3} & \text{if } \frac{1}{2} < H < \frac{3}{4},
	\end{cases}
	\end{align*}
	where we used the fact that if $H \neq \frac{1}{2}$, one has $\rho_H(r) \sim \frac{H(2H-1)}{a^2}\,r^{2H-2}$ for $r \to +\infty$.
\end{proof}

\begin{theorem}\label{thm:Q-cont}
	If $0 < H < \frac{3}{4}$ and $N \sim \mathcal{N}(0,1)$, then
	\[
	\dW\left(\frac{\sqrt{T}}{\sigma_H}\left(\bar Q_T - f_H(a)\right), N\right) \lesssim
	\begin{cases}
	T^{-1/2} & \text{if } 0 < H < 1/2 \\
	T^{2H-3/2} & \text{if } 1/2 \leqslant H < 3/4.
	\end{cases}
	\]
	where $\sigma_H^2 := 4\int_0^{+\infty}\rho_H^2(r)\,dr$, with $\rho_H(t) = \mathbb{E}[Z_t^H Z_0^H]$ defined above. Consequently, if $0 < H < 3/4$, then, as $T \to \infty$
	\[
	\sqrt{T}\left(\bar Q_T - f_H(a)\right) \xrightarrow{\text{law}} \mathcal{N}(0, \sigma_H^2).
	\]
\end{theorem}

\begin{proof}
	Using standard properties of the Wasserstein distance (see for instance \cite[Lemma~9]{EsSebaiyViens2019}) and Proposition 5, Lemma 6 and Proposition 2, we obtain
	\begin{align*}
	\dW\left(\frac{\sqrt{T}}{\sigma_H}\left(\bar Q_T - f_H(a)\right), N\right) \leqslant{}& C\left(\sqrt{T}\|R_H(T)\|_{L^1(\Omega)} + |v_T - \sigma_H^2|\right) \\
	&+ C
	\begin{cases}
	T^{-1/2} & \text{if } H \in \left(0,\frac{2}{3}\right) \\
	T^{-1/2}\log(T)^2 & \text{if } H = \frac{2}{3} \\
	T^{6H-9/2} & \text{if } H \in \left(\frac{2}{3},\frac{3}{4}\right).
	\end{cases}
	\end{align*}
	The desired result follows.
\end{proof}

\subsubsection{Berry--Ess\'een bounds for CLT for $\hat a_T$}

Notice that the estimator $\hat a_T$ can be written as $\hat a_T = f_H^{-1}(\bar Q_T)$, with $f_H^{-1}(x) = \left(\frac{H\Gamma(2H)}{x}\right)^{1/(2H)}$, besides by the mean-value theorem there exists a random variable $\zeta_T \in [|\bar Q_T, f_H(a)|]$, where the notation $[|a,b|]$ for two reals $a,b$ means $[a,b]$ if $a \leqslant b$ and $[b,a]$ if $b \leqslant a$, such that
\[
(\hat a_T - a) = (f_H^{-1})'(\zeta_T) \times \left(\bar Q_T - f_H(a)\right).
\]
Therefore
\begin{align*}
&\dW\left(\frac{\sqrt{T}}{g_H(a)\sigma_H}(\hat a_T - a), \mathcal{N}(0,1)\right) \\
&\quad\leqslant \dW\left(\frac{\sqrt{T}}{g_H(a)\sigma_H}(\hat a_T - a), \frac{\sqrt{T}}{\sigma_H}(\bar Q_T - f_H(a))\right) \\
&\qquad + \dW\left(\frac{\sqrt{T}}{\sigma_H}(\bar Q_T - f_H(a)), \mathcal{N}(0,1)\right)
\end{align*}
where $g_H(a) = (f_H^{-1})'(f_H(a)) = -\frac{a^{1+2H}}{2\Gamma(2H)H^2}$. Since $f_H^{-1}$ is twice continuously differentiable, applying the mean-value theorem again, there exists a random variable $\delta_T \in [|\zeta_T, f_H(a)|] \subset [|\bar Q_T, f_H(a)|]$.
\begin{align*}
&\dW\left(\frac{\sqrt{T}}{g_H(a)\sigma_H}(\hat a_T - a), \frac{\sqrt{T}}{\sigma_H}(\bar Q_T - f_H(a))\right) \times |g_H(a)\sigma_H| \\
&\quad \leqslant \sqrt{T}\,\mathbb{E}\left[\left|(\bar Q_T - f_H(a))\left((f_H^{-1})'(\zeta_T) - (f_H^{-1})'(f_H(a))\right)\right|\right] \\
&\quad \leqslant \sqrt{T}\,\mathbb{E}\left[|(\bar Q_T - f_H(a))| |(f_H^{-1})''(\delta_T)(\zeta_T - f_H(a))|\right] \\
&\quad \leqslant \sqrt{T}\,\mathbb{E}\left[|(\bar Q_T - f_H(a))|^2 |(f_H^{-1})''(\delta_T)|\right] \\
&\quad \leqslant \sqrt{T}\,\mathbb{E}\left[|(\bar Q_T - f_H(a))|^{2p}\right]^{1/p} \times \mathbb{E}\left[|(f_H^{-1})''(\delta_T)|^{p'}\right]^{1/p'},
\end{align*}
where we used H\"older's inequality with $p, p'$ two reals greater than 1 such that $\frac{1}{p}+\frac{1}{p'}=1$. Moreover, by the hypercontractivity property for multiple integrals, there exists a constant $C(p)$ such that
\begin{align}
\sqrt{T}\,\mathbb{E}\left[|(\bar Q_T - f_H(a))|^{2p}\right]^{1/p} &\leqslant C(p) \times \sqrt{T}\,\mathbb{E}\left[|(\bar Q_T - f_H(a))|^2\right] \notag\\
&\leqslant 2C(p)\sqrt{T}\,\mathbb{E}\left[|\bar Q_T - \mathbb{E}[\bar Q_T]|^2\right] + 2C(p)\sqrt{T}\left(\mathbb{E}[\bar Q_T] - f_H(a)\right)^2 \tag{12}
\end{align}

\textbf{The term $\mathbb{E}[|\bar Q_T - \mathbb{E}[\bar Q_T]|^2]$:} From the expression of the estimator $\bar Q_T$ defined in (4), we can write
\[
\bar Q_T - \mathbb{E}[\bar Q_T] = A(T) - B(T),
\]
where
\[
A(T) := \frac{1}{T}\int_0^T \left(X_t^2 - \mathbb{E}[X_t^2]\right) dt \quad \text{and} \quad B(T) := \left(\frac{1}{T}\int_0^T X_t dt\right)^2 - \mathbb{E}\left[\left(\frac{1}{T}\int_0^T X_t dt\right)^2\right].
\]
Let us now estimate $\mathbb{E}[A^2(T)]$. Using the fact that the process $X$ can be written as
\[
X_t = k(t) + G_t^H
\]
where $k(t) := b(1-e^{-at})$ and $G_t^H := \int_0^t e^{-a(t-s)} dB^H_s$, we get that
\[
A(T) = \frac{2}{T}\int_0^T k(t) G_t^H dt + \frac{1}{T}\int_0^T\left((G_t^H)^2 - \mathbb{E}[(G_t^H)^2]\right) dt.
\]
Hence, we have the estimate
\[
\mathbb{E}[A^2(T)] \leqslant 2 \times \mathbb{E}\left[\left(\frac{2}{T}\int_0^T k(t)G_t^H dt\right)^2\right] + 2 \times \mathbb{E}\left[\left(\frac{1}{T}\int_0^T\left((G_t^H)^2 - \mathbb{E}[(G_t^H)^2]\right)dt\right)^2\right].
\]
On the other hand,
\begin{align*}
\mathbb{E}\left[\left(\frac{2}{T}\int_0^T k(t)G_t^H dt\right)^2\right]
&= \frac{4}{T^2}\int_0^T\int_0^T k(t)k(s)\mathbb{E}[G_t^H G_s^H]\,dt\,ds \\
&= \frac{4}{T^2}\int_0^T\int_0^T k(t)k(s)\rho_H(t-s)\,dt\,ds \\
&\quad - \frac{8}{T^2}\int_0^T\int_0^T e^{-as}k(t)k(s)\rho_H(t)\,dt\,ds \\
&\quad + \frac{4f_H(a)}{T^2}\int_0^T\int_0^T e^{-a(t+s)}k(t)k(s)\,dt\,ds,
\end{align*}
where we used the fact that $G_t^H = Z_t^H - e^{-at}Z_0^H$, where $Z^H$ is the process defined in (7) with covariance function $\rho_H$. Moreover
\begin{itemize}
	\item
	\begin{align*}
	\frac{4}{T^2}\int_0^T\int_0^T k(t)k(s)\rho_H(t-s)\,dt\,ds &= \frac{4b^2}{T^2}\int_0^T\int_0^T (1-e^{-as})(1-e^{-at})\rho_H(t-s)\,dt\,ds \\
	&\leqslant \frac{4b^2}{T^2}\int_0^T\int_0^T \rho_H(t-s)\,dt\,ds = \frac{4b^2}{T}\int_0^T\left(1-\frac{u}{T}\right)\rho_H(u)\,du \\
	&\leqslant C
	\begin{cases}
	T^{-1} & \text{if } 0 < H \leqslant \frac{1}{2}, \\
	T^{2H-2} & \text{if } \frac{1}{2} < H < 1.
	\end{cases}
	\end{align*}
	\item
	\begin{align*}
	\frac{1}{T^2}\int_0^T\int_0^T e^{-as}k(t)k(s)\rho_H(t)\,dt\,ds &= \frac{1}{T^2}\int_0^T k(t)\rho_H(t)\,dt \times \int_0^T k(s)e^{-as}\,ds \\
	&\leqslant C
	\begin{cases}
	T^{-2} & \text{if } 0 < H \leqslant \frac{1}{2}, \\
	T^{2H-3} & \text{if } \frac{1}{2} < H < 1.
	\end{cases}
	\end{align*}
	\item $\left(\dfrac{1}{T^2}\displaystyle\int_0^T e^{-at}k(t)\,dt\right)^2 \leqslant C\dfrac{1}{T^2}$.
\end{itemize}

It follows that
\[
\mathbb{E}\left[\left(\frac{2}{T}\int_0^T k(t)G_t^H dt\right)^2\right] \leqslant C
\begin{cases}
T^{-1} & \text{if } 0 < H \leqslant \frac{1}{2}, \\
T^{2H-2} & \text{if } \frac{1}{2} < H < 1.
\end{cases}
\]

On the other hand, we have by the product formula
\[
\frac{1}{T}\int_0^T\left((G_t^H)^2 - \mathbb{E}[(G_t^H)^2]\right) dt = \frac{1}{T}\int_0^T \left(I_1^2(g_t) - \mathbb{E}[I_1^2(g_t)]\right) dt = I_2\left(\frac{1}{T}\int_0^T g_t^{\otimes 2} dt\right)
\]
with $g_t(\cdot) = e^{-a(t-\cdot)}\mathbf{1}_{[0,t]}$, consequently
\begin{align*}
\mathbb{E}\left[\left(\int_0^T \left((G_t^H)^2 - \mathbb{E}[(G_t^H)^2]\right) dt\right)^2\right]
&= \frac{2}{T^2}\int_0^T\int_0^T \langle g_t^{\otimes 2}, g_s^{\otimes 2}\rangle\,dt\,ds = \frac{2}{T^2}\int_0^T\int_0^T (\langle g_t, g_s\rangle)^2\,dt\,ds \\
&= \frac{2}{T^2}\int_0^T\int_0^T \mathbb{E}[G_t^H G_s^H]^2\,dt\,ds.
\end{align*}

Therefore, using the fact that $\mathbb{E}[G_t^H G_s^H] = \rho_H(t-s) - e^{-as}\rho_H(t) - e^{-at}\rho_H(s) + e^{-a(t+s)}f_H(a)$, we get
\begin{align*}
\mathbb{E}\left[\left(\int_0^T \left((G_t^H)^2 - \mathbb{E}[(G_t^H)^2]\right) dt\right)^2\right]
&= \frac{2}{T^2}\int_0^T\int_0^T \rho_H(t-s)\,dt\,ds - \frac{4}{T^2}\int_0^T\int_0^T e^{-as}\rho_H(t)\,dt\,ds \\
&\quad + \frac{2f_H(a)}{T^2}\int_0^T\int_0^T e^{-a(t+s)}\,dt\,ds.
\end{align*}
On the other hand,
\[
\frac{2}{T^2}\int_0^T\int_0^T \rho_H(t-s)\,dt\,ds \leqslant C
\begin{cases}
T^{-1} & \text{if } 0 < H \leqslant \frac{1}{2}, \\
T^{2H-2} & \text{if } \frac{1}{2} < H < 1,
\end{cases}
\]
moreover,
\[
\frac{4}{T^2}\int_0^T\int_0^T e^{-as}\rho_H(t)\,dt\,ds \leqslant C
\begin{cases}
T^{-2} & \text{if } 0 < H \leqslant \frac{1}{2}, \\
T^{2H-3} & \text{if } \frac{1}{2} < H < 1,
\end{cases}
\]
and
\[
\frac{2f_H(a)}{T^2}\int_0^T\int_0^T e^{-a(t+s)}\,dt\,ds \leqslant \frac{C}{T^2}.
\]

Finally, we get
\[
\mathbb{E}\left[\left(\int_0^T \left((G_t^H)^2 - \mathbb{E}[(G_t^H)^2]\right) dt\right)^2\right] \leqslant C
\begin{cases}
T^{-1} & \text{if } 0 < H \leqslant \frac{1}{2}, \\
T^{2H-2} & \text{if } \frac{1}{2} < H < 1.
\end{cases}
\]

Hence for the term $A(T)$, we have
\[
\mathbb{E}[A^2(T)] \leqslant C
\begin{cases}
T^{-1} & \text{if } 0 < H \leqslant \frac{1}{2}, \\
T^{2H-2} & \text{if } \frac{1}{2} < H < 1.
\end{cases}
\]

For the term $B(T)$, we have
\[
\frac{1}{T}\int_0^T X_t dt = \frac{1}{T}\int_0^T k(t)\,dt + \frac{1}{T}\int_0^T G_t^H\,dt.
\]
Thus,
\[
B(T) = \frac{2}{T^2}\int_0^T\int_0^T k(s)G_t^H\,dt\,ds + \frac{1}{T^2}\int_0^T\int_0^T \left(G_t^H G_s^H - \mathbb{E}[G_t^H G_s^H]\right)\,dt\,ds.
\]
We get therefore the following estimate:
\begin{equation}
\mathbb{E}[B^2(T)] \leqslant 2\mathbb{E}\left[\left(\frac{2}{T^2}\int_0^T k(s)\,ds \int_0^T G_t^H dt\right)^2\right] + 2\mathbb{E}\left[\left(\frac{1}{T^2}\int_0^T\int_0^T \left(G_t^H G_s^H - \mathbb{E}[G_t^H G_s^H]\right) dt\,ds\right)^2\right]. \tag{13}
\end{equation}
For the left hand term, we have
\[
\mathbb{E}\left[\left(\frac{2}{T^2}\int_0^T k(s)\,ds \int_0^T G_t^H\,dt\right)^2\right] = \frac{4}{T^4}\int_0^T\int_0^T k(v)k(u)\,dv\,du \int_0^T\int_0^T \mathbb{E}[G_t^H G_s^H]\,dt\,ds.
\]
On the other hand, we have
\begin{align*}
\int_0^T\int_0^T \mathbb{E}[G_t^H G_s^H]\,dt\,ds ={}& \int_0^T\int_0^T \rho_H(t-s)\,dt\,ds - 2\int_0^T \rho_H(t)\,dt \int_0^T e^{-as}\,ds \\
&+ f_H(a)\int_0^T\int_0^T e^{-a(t+s)}\,dt\,ds.
\end{align*}
Similarly to some previous calculations, we get
\[
\mathbb{E}\left[\left(\frac{2}{T^2}\int_0^T k(s)\,ds \int_0^T G_t^H\,dt\right)^2\right] \leqslant C
\begin{cases}
T^{-1} & \text{if } 0 < H \leqslant \frac{1}{2}, \\
T^{2H-2} & \text{if } \frac{1}{2} < H < 1.
\end{cases}
\]
On the other hand,
\begin{align}
D_T &:= \frac{1}{T}\int_0^T G_t^H\,dt = \frac{1}{T}\int_0^T \left(\int_0^t e^{-a(t-u)} dB^H_u\right) dt = \frac{1}{T}\int_0^T \left(\int_u^T e^{-a(t-u)} dt\right) dB^H_u \tag{14}\\
&= \frac{1}{aT}\left(B_T^H - G_T^H\right). \tag{15}
\end{align}

Moreover, since $G_T^H = Z_T - e^{-aT}Z_0$, then using the stationarity of the process $Z$, we get for all $H \in (0,1)$,
\begin{equation}
\mathbb{E}[(G_T^H)^2] = (1+e^{-2aT})f_H(a) - 2e^{-aT}\rho_H(T) \to \frac{H\Gamma(2H)}{a^{2H}}, \tag{16}
\end{equation}
as $T \to +\infty$. Consequently, using the previous calculus
\begin{equation}
D_T = \frac{1}{aT}\left(B_T^H - G_T^H\right) = \frac{B_T^H}{aT} + O\left(\frac{1}{T}\right). \tag{17}
\end{equation}
Thus, $D_T^2 = \frac{(B_T^H)^2}{a^2T^2} + O(T^{H-2})$, hence $\mathbb{E}[D_T^2] = O(T^{2H-2})$, where we used the self-similarity property of $B^H$. Using the hypercontractivity property, it follows that there exists a constant $C$, such that
\[
\mathrm{Var}(D_T^2) \leqslant C\left(\mathbb{E}[D_T^2]\right)^2 = O(T^{4H-4}).
\]
It follows that the right-hand term of inequality (13) equals
\[
\mathbb{E}\left[\left(\frac{1}{T^2}\int_0^T\int_0^T \left(G_t^H G_s^H - \mathbb{E}[G_t^H G_s^H]\right) dt\,ds\right)^2\right] = \mathrm{Var}(D_T^2).
\]
Thus, for the term $B(T)$, we have
\[
\mathbb{E}[B^2(T)] \leqslant C
\begin{cases}
T^{-1} & \text{if } 0 < H \leqslant \frac{1}{2}, \\
T^{2H-2} & \text{if } \frac{1}{2} < H < 1.
\end{cases}
\]
Finally, we can conclude that
\[
\sqrt{T}\,\mathbb{E}[|\bar Q_T - \mathbb{E}[\bar Q_T]|^2] \leqslant C
\begin{cases}
T^{-1/2} & \text{if } 0 < H \leqslant \frac{1}{2}, \\
T^{2H-3/2} & \text{if } \frac{1}{2} < H < 1.
\end{cases}
\]
On the other hand, we have by Lemma 3, we have
\[
\sqrt{T}|\mathbb{E}[\bar Q_T] - f_H(a)| \leqslant C
\begin{cases}
T^{-3/2} & \text{if } 0 < H \leqslant \frac{1}{2}, \\
T^{4H-7/2} & \text{if } \frac{1}{2} < H < 1.
\end{cases}
\]
It follows from inequality (12) that for any $p \geqslant 1$, there exists a constant $C(p) > 0$ such that
\[
\sqrt{T}\,\mathbb{E}\left[|\bar Q_T - f_H(a)|^{2p}\right]^{1/p} \leqslant C(p)
\begin{cases}
T^{-1/2} & \text{if } 0 < H \leqslant \frac{1}{2}, \\
T^{2H-3/2} & \text{if } \frac{1}{2} < H < 1.
\end{cases}
\]
It remains to show that $\mathbb{E}[|(f_H^{-1})''(\delta_T)|^{p'}]$ is finite for some $p' \geqslant 1$. From the monotonicity of $(f_H^{-1})''$ it's enough to show that $\mathbb{E}[|\bar Q_T|^{-p}] < +\infty$ for some $p \geqslant 1$.

\section{Berry--Ess\'een bounds for $\hat b_T$}\label{sec:b-cont}

For the estimator $\hat b_T$, we have from the definition (3):
\begin{align}
\hat b_T &= \frac{1}{T}\int_0^T b(1-e^{-at})\,dt + \frac{1}{T}\int_0^T \left(\int_0^t e^{-a(t-u)} dB^H_u\right) dt \notag\\
&= b - \frac{1}{T}\left(\frac{1-e^{-aT}}{a}\right) + D_T, \tag{18}
\end{align}
where $D_T$ is defined in (14). It follows that $\hat b_T - b = D_T - \frac{1}{T}\left(\frac{1-e^{-aT}}{a}\right)$. Moreover, by (17), we have $D_T = \frac{1}{aT}(B_T^H - G_T^H)$, thus $\mathbb{E}[D_T^2] = \frac{\mathbb{E}[(B_T^H)^2]}{a^2T^2} + O(T^{H-2}) = \frac{T^{2H-2}}{a^2} + O(T^{H-2})$.

Therefore, the following Lemma holds.

\begin{lemma}
	Consider the random variable $D_T$ defined in (14), then for all $H \in (0,1)$, we have
	\[
	\left|T^{2-2H}\mathbb{E}[D_T^2] - \frac{1}{a^2}\right| \leqslant \frac{C}{T^{2-H}}.
	\]
	Consequently, we have $\displaystyle\lim_{T\to+\infty} T^{2-2H}\mathbb{E}[D_T^2] = \frac{1}{a^2}$.
\end{lemma}

Concerning the rate of convergence, we have the following estimate for the Wasserstein distance for the estimator $\hat b_T$.

\begin{theorem}\label{thm:b-cont}
	Consider the random variable $D_T$ defined in (14), then for all $H \in (0,1)$, we have:
	\[
	\dW\left(T^{1-H}(\hat b_T - b), \mathcal{N}\left(0, \frac{1}{a^2}\right)\right) = O\left(\frac{1}{T^H}\right).
	\]
	Consequently, the following CLT follows: for all $0 < H < 1$ and as $T \to \infty$
	\[
	T^{1-H}(\hat b_T - b) \xrightarrow{\text{law}} \mathcal{N}\left(0, \frac{1}{a^2}\right).
	\]
\end{theorem}

\begin{proof}
	In fact, for the convergence in law of the term $D_T$, we have from (17),
	\begin{align*}
	\dW\left(aT^{1-H}D_T, \mathcal{N}(0,1)\right) &= \dW\left(aT^{1-H}\left(\frac{1}{aT}(B_T^H + O(\tfrac{1}{T}))\right), \mathcal{N}(0,1)\right) \\
	&\leqslant \dW\left(T^{-H}B_T^H, \mathcal{N}(0,1)\right) + O(T^{-H}) = O(T^{-H}),
	\end{align*}
	since $B_T^H \stackrel{\text{law}}{=} T^H B_1^H$. It's easy to conclude, since from (18):
	\begin{align*}
	\dW\left(T^{1-H}(\hat b_T - b), \mathcal{N}\left(0, \frac{1}{a^2}\right)\right) &\leqslant O\left(\frac{1}{T^H}\right) + \dW\left(aT^{1-H}D_T, \mathcal{N}(0,1)\right) \\
	&= O\left(\frac{1}{T^H}\right).
	\end{align*}
\end{proof}
\section{Discrete-time observations}\label{sec:discrete}

We now observe the fractional Vasicek process $X$ of (1) at the instants $t_i=i\Delta_n$, $i=0,\dots,n$, with
$\Delta_n\to0$ and $T_n=n\Delta_n\to\infty$. Throughout we assume without loss of generality that $\Delta_n\le 1\le T_n$. The constant $C$ may change from line to line; it depends on $(a,b,H)$ (and on $p$ when moments of order $p$ are involved) but never on $n$.\\
The estimators are
\[
\hat b_n=\frac1n\sum_{i=0}^{n-1}X_{t_i},\qquad
\bQ_n=\frac1n\sum_{i=0}^{n-1}X_{t_i}^2-\hat b_n^{\,2},\qquad
\hat a_n=\Bigl(\frac{\bQ_n}{H\Gamma(2H)}\Bigr)^{-1/(2H)} .
\]
Consider the two discretization errors :
\[
E_1(n):=\hat b_n-\hat b_{T_n},\qquad E_2(n):=\bQ_n-\bQ(T_n),
\]
where $\hat b_{T_n}$ and $\bQ(T_n)$ are the continuous-time estimators of Section~\ref{sec:continuous} evaluated at $T_n=T_n$.

\paragraph{Rates.} Put
\begin{align*}
& \gamma_H :=\min\{H+\tfrac12,\,1\},\\
& \tau_H(T_n) :=\begin{cases}T_n^{-1/2}&0<H\le\frac12,\\[2pt] T_n^{H-1}&\frac12<H<1,\end{cases}\\
 & \nu_H(T_n) :=\sqrt{T_n}\,\tau_H(T_n)^2=\begin{cases}T_n^{-1/2}&0<H<\frac12,\\[2pt] T_n^{2H-3/2}&\frac12\le H<1.\end{cases}
\end{align*}
\paragraph{Covariance estimates.} The covariance $\rho_H(u)=\E[Z_u^HZ_0^H]$ of the stationary process $Z^H$ is controlled by Lemma~\ref{lem:rho} of Appendix~\ref{app:cov}; it is used below with $\varphi=\rho_H$ and $\varphi=\rho_H^2$, and $f_H(a)=\rho_H(0)$.

\subsection{The block identity and the linear error}

Since $X_t-X_s=a\int_s^t(b-X_u)du+B^H_t-B^H_s$ pathwise, integrating over a block $[t_i,t_{i+1}]$ gives the block identity :
\begin{equation}\label{eq:block}
\int_{t_i}^{t_{i+1}}(X_t-X_{t_i})\,dt=a\int_{t_i}^{t_{i+1}}w(u)\,(b-X_u)\,du+\xi_i,
\qquad w(u):=t_{i+1}-u,\quad \xi_i:=\int_{t_i}^{t_{i+1}}(B^H_t-B^H_{t_i})\,dt .
\end{equation}
For $H\ne\frac12$ the variables $\xi_i$ are not independent, but their dependence is controlled by the following lemma.

\begin{lemma}\label{lem:zeta}
\begin{enumerate}
\item[(a)] $\E\xi_i^2\le\Delta_n^{2H+2}$ for every $i$;
\item[(b)] for $|i-j|\ge2$,
\begin{equation}\label{eq:xi-cov}
|\E[\xi_i\xi_j]|\le H|2H-1|\,\Delta_n^{2H+2}\,(|i-j|-1)^{2H-2};
\end{equation}
\item[(c)] $\Var\bigl(\sum_{i=0}^{n-1}\xi_i\bigr)\le C\,n\,\Delta_n^{2H+2}\max\{1,n^{2H-1}\}$, and consequently :
\[
\Bigl\lVert\frac1{T_n}\sum_{i=0}^{n-1}\xi_i\Bigr\rVert_{L^2}\le C\,\Delta_n^{\gamma_H}\tau_H(T_n).
\]
\end{enumerate}
\end{lemma}

\begin{proof}
\begin{itemize}
\item (a) By Cauchy--Schwarz, $\E\xi_i^2\le\Delta_n\int_{t_i}^{t_{i+1}}\E(B^H_t-B^H_{t_i})^2dt=\Delta_n\int_0^{\Delta_n} s^{2H}ds \le \frac{{\Delta_n}^{2H+2}}{2H+1}$.
\\
\item (b) Since $\partial_x\partial_y\frac12\bigl(x^{2H}+y^{2H}-|x-y|^{2H}\bigr)=H(2H-1)|x-y|^{2H-2}$ off the diagonal, for disjoint intervals $[u,t]$ and $[v,s]$, we can write for all $H \in (0,1)$ :
\[
\E\bigl[(B^H_t-B^H_u)(B^H_s-B^H_v)\bigr]=H(2H-1)\int_u^t\!\!\int_v^s|x-y|^{2H-2}dy\,dx .
\]

If $|i-j| \geq 2$, we apply this to $u=t_i$ and $v= t_j$, then 
$$\mathbb{E}[\xi_i \xi_j] = \int_{t_{i}}^{t_{j+1}} \int_{t_j}^{t_{j+1}} \mathbb{E}[(B_t^{H} - B_{t_i}^{H}) (B_s^{H} - B_{t_j}^{H})] ds dt,$$
\begin{itemize}
\item[\textbullet] If $|i-j| \geq 2$, the two blocs are separated by at least one block. That is for $x \in [t_i,t_{i+1}]$ and $y \in [t_j, t_{j+1}]$, we have $|x-y| \geq (|i-j|-1) \Delta_n$.
\item[\textbullet] 	As $2H-2 < 0$, the function $r \mapsto r^{2H-2}$ is decreasing, that is $|x-y|^{2H-2} \leq ((|i-j|-1) \Delta_n)^{2H-2}$ for all $H \in (0,1)$.
That is the following estimate follows : for all $i,j$ such that  $|i-j| \geq 2$
$$|\mathbb{E}[\xi_i \xi_j]| \leq H|2H-1| \Delta_n^{2H+2} (|i-j|-1)^{2H-2}.$$
\end{itemize}

\item (c) We have
\[
\Var\Bigl(\sum_i\xi_i\Bigr)=\sum_{i,j}\E[\xi_i\xi_j]=\sum_{\substack{i,j\\|i-j|\le1}}\E[\xi_i\xi_j]+\sum_{\substack{i,j\\|i-j|\ge2}}\E[\xi_i\xi_j].
\]
\emph{First sum.}  By Cauchy--Schwarz and (a), $|\E[\xi_i\xi_j]|\le\Delta_n^{2H+2}$, so
\[
\sum_{|i-j|\le1}|\E[\xi_i\xi_j]|\le3n\Delta_n^{2H+2}.
\]
\emph{Second sum.} By \eqref{eq:xi-cov},

\begin{align*}
\sum_{|i-j|\ge2}|\E[\xi_i\xi_j]|
&\le H|2H-1|\,\Delta_n^{2H+2}\sum_{|i-j|\ge2}(|i-j|-1)^{2H-2}\\
&=2H|2H-1|\,\Delta_n^{2H+2}\sum_{i=0}^{n-3}\ \sum_{j=i+2}^{n-1}(j-i-1)^{2H-2}.
\end{align*}
Hence
\[
\sum_{|i-j|\ge2}|\E[\xi_i\xi_j]|\le2H|2H-1|\,\Delta_n^{2H+2}\sum_{r=1}^{n-2}\ \sum_{i=0}^{n-r-2}r^{2H-2}
=2H|2H-1|\,\Delta_n^{2H+2}\sum_{r=1}^{n-2}(n-r-1)\,r^{2H-2}.
\]
Using a comparison between series and integrals, we get, according to the value of $H$:
\begin{itemize}
\item if $H<\frac12$:  $\displaystyle\sum_{|i-j|\ge2}|\E[\xi_i\xi_j]|\le4H(1-H)\,n\,\Delta_n^{2H+2}$;
\item if $H>\frac12$:  $\displaystyle\sum_{|i-j|\ge2}|\E[\xi_i\xi_j]|\le2H\,n^{2H}\Delta_n^{2H+2}$;
\end{itemize}
Altogether
\[
\Var\Bigl(\sum_i\xi_i\Bigr)\le C\,n\,\Delta_n^{2H+2}\max\{1,n^{2H-1}\}.
\]
Now we divide by $T_n=n\Delta_n$.
For $H\le\frac12$,
\[
\frac1T_n\Bigl\lVert\sum_i\xi_i\Bigr\rVert_{L^2}\le C\,\frac1{n\Delta_n}\,\sqrt n\,\Delta_n^{H+1}=C\,\frac{\Delta_n^{H}}{\sqrt n}=C\,\Delta_n^{H+1/2}\,T_n^{-1/2}.
\]
For $H>\frac12$,
\[
\frac1{T_n}\Bigl\lVert\sum_i\xi_i\Bigr\rVert_{L^2}\le C\,\frac{n^{H}\Delta_n^{H+1}}{n\Delta_n}=C\,n^{H-1}\Delta_n^{H}=C\,\Delta_n\,T_n^{H-1}.
\]

The desired results follows by definition of the rate $\gamma_H$ and $\tau_H$.
\end{itemize}
\end{proof}

\begin{lemma}[Linear error]\label{lem:linear}
For every $n$, $\ \norm{E_1(n)}{L^2}\le C\,\Delta_n^{\gamma_H}\,\tau_H(T_n)$, and $E_1(n)$ is Gaussian, so $\norm{E_1(n)}{L^p}\le C_p\norm{E_1(n)}{L^2}$ for all $p\ge1$.
\end{lemma}

\begin{proof}
Since $X=b+Z^H+Y$, we can write
\[
E_1(n)=(\bar Z_n^H-\bar Z_{T_n}^H)+(\bar Y_n-\bar Y_{T_n}),
\]
where $\bar Z_n^H=\frac1n\sum_{i<n}Z_{t_i}^H$, $\bar Z_{T_n}^H=\frac1T_n\int_0^TZ_t^Hdt$, and similarly $\bar Y_n,\bar Y_{T_n}$. Recall that $E_1(n)=\hat b_n-\hat b_{T_n}=\frac1n\sum_{i=0}^{n-1}X_{t_i}-\frac1{T_n}\int_0^{T_n}X_t\,dt$.

\medskip
\noindent\textbf{(1) The $Y$-part.} We have $\bar Y_n-\bar Y_{T_n}=-(b+Z_0^H)(e_n-e_{T_n})$, with  $e_n=\frac1n\sum_{i=0}^{n-1}e^{-at_i}$ and $e_{T_n}=\frac1T_n\int_0^{T_n}e^{-at}dt$. On the other hand,
\[
e_n-e_{T_n}=\frac1n\sum_{i=0}^{n-1}e^{-at_i}-\frac1{T_n}\int_0^Te^{-at}dt=\frac1{n\Delta_n}\sum_{i=0}^{n-1}\int_{t_i}^{t_{i+1}}\bigl(e^{-at_i}-e^{-at}\bigr)dt .
\]
Since $t\mapsto-e^{-at}$ is increasing, we have
\begin{align*}
|e_n-e_{T_n}|&\le\frac1{n\Delta_n}\sum_{i=0}^{n-1}\int_{t_i}^{t_{i+1}}\bigl(e^{-at_i}-e^{-at_{i+1}}\bigr)dt
=\frac{\Delta_n}{n\Delta_n}\sum_{i=0}^{n-1}e^{-at_i}\bigl(1-e^{-a\Delta_n}\bigr)\\
&=\frac{\Delta_n\,(1-e^{-a\Delta_n})}{n\Delta_n}\cdot\frac{1-e^{-aT_n}}{1-e^{-a\Delta_n}}
=\frac{1-e^{-aT_n}}{n}\le\frac{\Delta_n}{T_n}\le\Delta_n^{\gamma_H}\tau_H(T_n),
\end{align*}
 Hence
\[
\norm{\bar Y_n-\bar Y_{T_n}}{L^2}\le\norm{b+Z_0^H}{L^2}\,\Delta_n^{\gamma_H}\tau_H(T_n)=\bigl(b^2+f_H(a)\bigr)^{1/2}\,\Delta_n^{\gamma_H}\tau_H(T_n).
\]

\medskip
\noindent\textbf{(2) The $Z^H$-part.} Pathwise (as for the proof of \eqref{eq:block}), for $t\in[t_i,t_{i+1}]$,
\[
Z_t^H-Z_{t_i}^H=-a\int_{t_i}^tZ_u^H\,du+\bigl(B^H_t-B^H_{t_i}\bigr).
\]
That is,
\begin{align*}
\bar Z_n^H-\bar Z_{T_n}^H&=\frac1{T_n}\sum_{i=0}^{n-1}\int_{t_i}^{t_{i+1}}(Z_{t_i}^H-Z_t^H)\,dt
=\frac1{T_n}\sum_{i=0}^{n-1}\int_{t_i}^{t_{i+1}}\Bigl(a\int_{t_i}^tZ_u^H\,du-(B^H_t-B^H_{t_i})\Bigr)dt\\
&{=}\frac aT\sum_{i=0}^{n-1}\int_{t_i}^{t_{i+1}}(t_{i+1}-u)\,Z_u^H\,du-\frac1{T_n}\sum_{i=0}^{n-1}\xi_i,
\end{align*}
with $\xi_i=\int_{t_i}^{t_{i+1}}(B^H_t-B^H_{t_i})\,dt$. Set $w_i(u)=t_{i+1}-u$ for $t_i\le u\le t_{i+1}$; then
\[
\bar Z_n^H-\bar Z_{T_n}^H=\frac aT\sum_{i=0}^{n-1}\int_{t_i}^{t_{i+1}}w_i(u)\,Z_u^H\,du-\frac1{T_n}\sum_{i=0}^{n-1}\xi_i .
\]
On the other hand, using Lemma~\ref{lem:rho}(a) along with the fact that $|w_i|\le\Delta_n$ for all $i=0,\dots,n-1$, we get
\[
\frac aT\Bigl\lVert\sum_{i=0}^{n-1}\int_{t_i}^{t_{i+1}}w_i(u)Z_u^H\,du\Bigr\rVert_{L^2}\le C\,\Delta_n\,T_n^{-1/2}\max\{1,T_n^{H-1/2}\}=C\,\Delta_n\,\tau_H(T_n)\le C\,\Delta_n^{\gamma_H}\tau_H(T_n).
\]
We already proved (Lemma~\ref{lem:zeta}(c)) that
\[
\Bigl\lVert\frac1{T_n}\sum_{i=0}^{n-1}\xi_i\Bigr\rVert_{L^2}\le C\,\Delta_n^{\gamma_H}\tau_H(T_n).
\]
Consequently, $\norm{\bar Z_n^H-\bar Z_{T_n}^H}{L^2}\le C\,\Delta_n^{\gamma_H}\tau_H(T_n)$.

\medskip
\noindent\textbf{Conclusion.} By Minkowski's inequality,
\[
\norm{E_1(n)}{L^2}\le\norm{\bar Y_n-\bar Y_{T_n}}{L^2}+\norm{\bar Z_n^H-\bar Z_{T_n}^H}{L^2}\le C\,\Delta_n^{\gamma_H}\tau_H(T_{n}).
\]
Moreover $E_1(n)$ is Gaussian so by hypercontractivity \eqref{eq:hypercontractivity} $\norm{E_1(n)}{L^p}\le C_p\norm{E_1(n)}{L^2}$ for every $p\ge1$.
\end{proof}

\subsection{The quadratic error}

Let $S_n:=\frac1n\sum_{i=0}^{n-1}(Z_{t_i}^H)^2-\frac1{T_n}\int_0^T(Z_t^H)^2dt$.

\begin{lemma}[Riemann sums of $(Z^H)^2$]\label{lem:S}
$S_n \in \mathcal{H}^{B^H}_2$, $\E S_n=0$, and $\norm{S_n}{L^2}\le C\,\Delta_n^{\gamma_H}\,T_n^{-1/2}$.
\end{lemma}

\begin{proof}
Since $\frac1n=\frac{\Delta_n}{T_n}$, we have
\begin{align*}
\frac1n\sum_i(Z_{t_i}^H)^2&=\frac1{T_n}\sum_i\int_{t_i}^{t_{i+1}}(Z_{t_i}^H)^2dt,\\
\frac1{T_n}\int_0^{T_n}(Z_t^H)^2dt&=\frac1{T_n}\sum_i\int_{t_i}^{t_{i+1}}(Z_t^H)^2dt .
\end{align*}
Hence
\[
S_n=\frac1{T_n}\sum_{i=0}^{n-1}\int_{t_i}^{t_{i+1}}\bigl((Z_{t_i}^H)^2-(Z_t^H)^2\bigr)\,dt .
\]
The process $Z^H$ is centered Gaussian, so $Z_t^H=I_1(\varepsilon_t)$ with $\langle\varepsilon_u,\varepsilon_v\rangle=\rho_H(u-v)$.\\
By the product formula \eqref{eq:product11}, $I_1(g)^2=I_2(g\otimes g)+\|g\|^2$, hence
\[
(Z_t^H)^2=I_2(\varepsilon_t^{\otimes2})+f_H\qquad(\text{since }\|\varepsilon_t\|^2=\rho_H(0)=f_H).
\]
The constants $f_H$ cancel in $S_n$, so
\[
S_n=I_2(h),\qquad h=\frac1{T_n}\sum_{i=0}^{n-1}\int_{t_i}^{t_{i+1}}\bigl(\varepsilon_{t_i}^{\otimes2}-\varepsilon_t^{\otimes2}\bigr)\,dt .
\]
In particular $S_n\in\mathcal H^{B^H}_2$ and $\E S_n=0$.\\
By the isometry property \eqref{eq:isometry}, $\E[I_2(h)^2]=2\|h\|^2$.
Expanding $\|h\|^2$ with
$\langle\varepsilon_u^{\otimes2},\varepsilon_v^{\otimes2}\rangle=\langle\varepsilon_u,\varepsilon_v\rangle^2=\rho_H(u-v)^2$, we get :
\[
\E S_n^2=\frac2{T_n^2}\sum_{i,j=0}^{n-1}\int_{t_i}^{t_{i+1}}\!\!\int_{t_j}^{t_{j+1}}D_{ij}(t,s)\,ds\,dt,
\]
\[
D_{ij}(t,s)=\rho_H(t_i-t_j)^2-\rho_H(t_i-s)^2-\rho_H(t-t_j)^2+\rho_H(t-s)^2.
\]

\textbf{Near the diagonal ($|i-j|\le1$).}
\begin{itemize}
\item
\[
D_{ij}=[\rho_H(t_i-t_j)^2-\rho_H(0)^2]-[\rho_H(t_i-s)^2-\rho_H(0)^2]-[\rho_H(t-t_j)^2-\rho_H(0)^2]+[\rho_H(t-s)^2-\rho_H(0)^2] .
\]
\item When $|i-j|\le1$, all four arguments have absolute value $\le2\Delta_n$. Since $\rho_H$ is even,
Lemma~\ref{lem:rho}(b) gives $|\rho_H(u)^2-\rho_H(0)^2|\le C|u|^{2H}$, hence $|D_{ij}|\le C\Delta_n^{2H}$.
\item There are at most $3n$ such pairs, each integrated over an area $\Delta_n^2$, so their contribution is at most
\[
\frac2{T_n^2}\cdot3n\cdot\Delta_n^2\cdot C\Delta_n^{2H}=C\,\frac{n\Delta_n^{2H+2}}{T_n^2}=C\,\frac{\Delta_n^{2H+1}}{T_n}.
\]
\end{itemize}

\item \textbf{Far from the diagonal ($k=|i-j|\ge2$), case $H\ne\frac12$.}
\begin{itemize}
\item $D_{ij}$ is the rectangular increment of $(x,y)\mapsto\rho_H(x-y)^2$ over $[t_i,t]\times[t_j,s]$.
Since $\partial_x\partial_y\rho_H(x-y)^2=-(\rho_H^2)''(x-y)$,
\[
D_{ij}(t,s)=-\int_{t_i}^t\!\!\int_{t_j}^s(\rho_H^2)''(x-y)\,dy\,dx .
\]
\item For $k\ge2$ the rectangle never meets the diagonal: $|x-y|\ge(k-1)\Delta_n>0$, where $\rho_H$ is smooth.
By Lemma~\ref{lem:rho}, $|(\rho_H^2)''|\le C\,m$, and $m$ is non-increasing; the integration area is $\le\Delta_n^2$, so
\[
|D_{ij}(t,s)|\le C\,\Delta_n^2\,m\bigl((k-1)\Delta_n\bigr).
\]
\item Summing over the pairs, we have
\[
\frac2{T_n^2}\sum_{|i-j|\ge2}\iint|D_{ij}|
\le\frac C{T_n^2}\,n\,\Delta_n^4\sum_{k\ge2}m\bigl((k-1)\Delta_n\bigr)
=\frac C{T_n}\,\Delta_n^3\sum_{r\ge1}m(r\Delta_n)\qquad(r=k-1).
\]
\item Let $J_n:=\lfloor 1/\Delta_n\rfloor$; since $\Delta_n\le1$, we have $J_n\ge1$ and $J_n\ge\frac1{2\Delta_n}$. Splitting the sum at $r=J_n$ according to the definition of $m$ (note that $r\Delta_n\le1$ iff $r\le J_n$), we have:
\begin{align*}
& \frac C{T_n}\Bigl(\Delta_n^{3}\Delta_n^{2H-2}\sum_{r\le J_n}r^{2H-2}+\Delta_n^{3}\Delta_n^{4H-6}\sum_{r>J_n}r^{4H-6}\Bigr) \\
& =\frac C{T_n}\Bigl(\Delta_n^{2H+1}\sum_{r\le J_n}r^{2H-2}+\Delta_n^{4H-3}\sum_{r>J_n}r^{4H-6}\Bigr).
\end{align*}
\item \emph{The tail $r>J_n$.} Since $4H-6<-1$, we have $\sum_{r>J_n}r^{4H-6}\le\int_{J_n}^\infty x^{4H-6}dx=\frac{J_n^{4H-5}}{5-4H}\le C\Delta_n^{5-4H}$.
The second term is therefore $\le C\,\Delta_n^{4H-3}\Delta_n^{5-4H}/T_n=C\,\Delta_n^2/T_n$.
\item \emph{The part $r\le J_n$.} If $H<\frac12$, the series $\sum r^{2H-2}$ converges, so the term is $\le C\Delta_n^{2H+1}/T_n$.
If $H>\frac12$,
\[
\sum_{r\le J_n}r^{2H-2}\le\int_0^{J_n}x^{2H-2}dx\le C\Delta_n^{1-2H},
\]
so the term is $\le C\Delta_n^{2H+1}\Delta_n^{1-2H}/T_n=C\Delta_n^2/T_n$.
\end{itemize}

\item \textbf{Conclusion.}
\begin{itemize}
\item Adding the contributions of the near-diagonal and far-from-diagonal pairs, we get $\E S_n^2\le C\,\Delta_n^{2H+1}/T_n$ if $H<\frac12$
and $\E S_n^2\le C\,\Delta_n^{2}/T_n$ if $H\ge\frac12$ (for $H>\frac12$, $\Delta_n^{2H+1}\le\Delta_n^2$ since $\Delta_n\le1$).
\item For $H<\frac12$ we also have $\Delta_n^2\le\Delta_n^{2H+1}$. In all cases
\[
\E S_n^2\le C\,\frac{\Delta_n^{\min\{2H+1,\,2\}}}{T_n}=C\,\frac{\Delta_n^{2\gamma_H}}{T_n}.
\]
The desired result follows.
\end{itemize}

\end{proof}

\begin{lemma}[Quadratic error]\label{lem:quad}
For every $p\ge1$,
\[
\norm{E_2(n)}{L^p}\le C_p\Bigl(T_n^{-1}+\Delta_n^{\gamma_H}\bigl(T_n^{-1/2}+\tau_H(T_n)^2\bigr)\Bigr).
\]
In particular, if $H<\frac34$, $\ \norm{E_2(n)}{L^p}\le C_p\,(T_n^{-1}+\Delta_n^{\gamma_H}T_n^{-1/2})$.
\end{lemma}

\begin{proof}
Set $\Lambda_t:=X_t-b=Z_t^H+Y_t$. Subtracting $b$ we get :
\[
E_2(n)=\Bigl[\frac1n\sum_i\Lambda_{t_i}^2-\frac1{T_n}\!\int_0^{T_n}\!\Lambda_t^2dt\Bigr]-\bigl(\bar\Lambda_n^{\,2}-\bar\Lambda_{T_n}^{\,2}\bigr).
\]
Expanding $\Lambda^2=(Z^H)^2+2Z^HY+Y^2$, the first bracket equals $S_n$ plus
$2\bigl[\frac1n\sum_iZ_{t_i}^HY_{t_i}-\frac1{T_n}\int Z^HY\bigr]+\bigl[\frac1n\sum_iY_{t_i}^2-\frac1{T_n}\int Y^2\bigr]$.
By Hölder's inequality, $\|Z_t^HY_t\|_p+\|Y_t^2\|_p\le Ce^{-at}$, and $\frac1n\sum_{i<n}e^{-at_i}\le\frac{1+a}{aT_n}$, $\frac1{T_n}\int_0^Te^{-at}dt\le\frac1{aT}$ (since $\Delta_n\le1$); thus the last two groups are $O(T_n^{-1})$ in $L^p$. Lemma~\ref{lem:S} and hypercontractivity \eqref{eq:hypercontractivity}  give $\norm{S_n}{p}\le C\Delta_n^{\gamma_H}T_n^{-1/2}$.

For the second part, $\bar\Lambda_n^2-\bar\Lambda_{T_n}^2=(\bar\Lambda_n-\bar\Lambda_{T_n})(\bar\Lambda_n+\bar\Lambda_{T_n})$ is a product of two Gaussian variables. By Lemma~\ref{lem:linear} and the $Y$-estimate there, $\norm{\bar\Lambda_n-\bar\Lambda_{T_n}}{2p}\le C\Delta_n^{\gamma_H}\tau_H(T_n)$, while $\norm{\bar Z_{T_n}^H}{2}^2=T_n^{-2}\iint\rho_H(u-v)\le C\tau_H(T_n)^2$ gives $\norm{\bar\Lambda_n+\bar\Lambda_{T_n}}{2p}\le C\tau_H(T_n)$ (using $\bar Z_n^H=\bar Z_{T_n}^H+(\bar Z_n^H-\bar Z_{T_n}^H)$). So this term is $\le C\Delta_n^{\gamma_H}\tau_H(T_n)^2$. Finally $\tau_H(T_n)^2\le T_n^{-1/2}$ when $H<\frac34$.
\end{proof}


\subsection{Consistency}

\begin{theorem}\label{thm:consistency-discrete}
Let $H\in(0,1)$, $\Delta_n\to0$ and $T_n\to\infty$.
\begin{enumerate}
\item[(a)] $(\hat a_n,\hat b_n)\to(a,b)$ in probability.
\item[(b)] If moreover $T_n\ge n^{\delta}$ for some $\delta>0$ (for instance $\Delta_n=n^{-\alpha}$, $0<\alpha<1$), then $(\hat a_n,\hat b_n)\to(a,b)$ almost surely.
\end{enumerate}
\end{theorem}

\begin{proof}
By Theorem~\ref{thm:consistency}, $\hat b_{T_n}\to b$ and $\bQ(T_n)\to f_H(a)$ a.s. as $T_n\to\infty$.
\begin{itemize}
\item \textbf{(a) Convergence in probability.}

By Lemmas~\ref{lem:linear} and~\ref{lem:quad}, $E_1(n)\to0$ and $E_2(n)\to0$ in $L^2$, hence in probability, so $\hat b_n\to b$ and $\bQ_n\to f_H(a)$ in probability.\\
Moreover $\bQ_n>0$ a.s.: for $n\ge2$, $\bQ_n=0$ forces $X_{t_1}=X_{t_0}=0$, and $X_{t_1}=b(1-e^{-a\Delta_n})+G_{t_1}$ is a non-degenerate Gaussian variable.
Finally, since the function $x\mapsto(x/(H\Gamma(2H)))^{-1/(2H)}$ is continuous on $(0,\infty)$, $\hat a_n\to a$ in probability.

\item \textbf{(b) Almost sure convergence.}
Both $E_1(n)$ and $E_2(n)$ live in a finite sum of chaoses, so by hypercontractivity \eqref{eq:hypercontractivity} all their $L^p$ norms are controlled by Lemmas~\ref{lem:linear}--\ref{lem:quad}: $\norm{E_1}{p}+\norm{E_2}{p}\le C_p\,T_n^{-\kappa}$ with $\kappa=\min\{\frac12,1-H\}>0$ (recall $\Delta_n\le1$ and $\tau_H(T_n)\le T_n^{-\kappa}$).\\
By Markov's inequality, for every $\varepsilon>0$,
\[
\Prob(|E_i(n)|>\varepsilon)\le C_p\,\varepsilon^{-p}\,n^{-\delta\kappa p},
\]
which is summable for $p>1/(\delta\kappa)$. Borel--Cantelli's Lemma gives $E_i(n)\to0$ a.s., which proves (b).
\end{itemize}
\end{proof}
\subsection{Berry--Ess\'een bounds}
The Berry Esséen bound for the estimator $\hat a_n$ requires uniform negative moments of $\bQ_n$ which is the next proposition, its proof can be found in Appendix.

\begin{proposition}[Negative moments, discrete case]\label{prop:neg-discrete}
Let $\Delta_n\le1$ and $T_n\to\infty$. For every $p\ge1$ there is $n_0=n_0(a,H,p)$ such that $\sup_{n\ge n_0}\E[\bQ_n^{-p}]<\infty$. 
\end{proposition}
We use the elementary coupling inequality: for real random variables $U,V$,
\begin{equation}\label{eq:coupling}
\dW(\mathrm{law}(U+V),\mathrm{law}(U))\le\E|V| .
\end{equation}
\begin{theorem}[Berry--Ess\'een bounds, discrete case]\label{thm:BE-discrete}
Let $N\sim\mathcal N(0,1)$ and let $\sigma_H^2=4\int_0^\infty\rho_H^2$, $g_H(a)=(f_H^{-1})'(f_H(a))$ be as in Section~\ref{sec:continuous}.
\begin{enumerate}
\item[(i)] If $0<H<\frac34$, there are $C$ and $n_0$ such that for all $n\ge n_0$
\begin{gather*}
\dW\Bigl(\frac{\sqrt{T_n}}{\sigma_H}\bigl(\bQ_n-f_H(a)\bigr),N\Bigr)\le C\,\mathcal R_n,\\
\dW\Bigl(\frac{\sqrt{T_n}}{g_H(a)\sigma_H}(\hat a_n-a),N\Bigr)\le C\,\mathcal R_n,
\end{gather*}
where $\mathcal R_n:=\nu_H(T_n)+\Delta_n^{\gamma_H}$, that is
\[
\mathcal R_n=\left\lbrace
\begin{array}{ll}
T_n^{-1/2}+\Delta_n^{H+1/2} & \text{if } 0<H<\frac12,\\[4pt]
T_n^{2H-3/2}+\Delta_n & \text{if } \frac12\le H<\frac34.
\end{array}\right.
\]
\item[(ii)] If $0<H<1$,
\[
\dW\Bigl(T_n^{1-H}(\hat b_n-b),\mathcal N\bigl(0,\tfrac1{a^2}\bigr)\Bigr)\le
C\,\left\lbrace
\begin{array}{ll}
T_n^{-H}+T_n^{1/2-H}\Delta_n^{H+1/2} & \text{if } 0<H\le\frac12,\\[4pt]
T_n^{-H}+\Delta_n & \text{if } \frac12<H<1,
\end{array}\right.
\]
that is, $\dW\le C\bigl(T_n^{-H}+T_n^{1-H}\Delta_n^{\gamma_H}\tau_H(T_n)\bigr)$.
\end{enumerate}
\end{theorem}

\begin{proof}
(i) \emph{The estimator $\bQ_n$.} Apply \eqref{eq:coupling} with $U=\frac{\sqrt{T_n}}{\sigma_H}(\bQ(T_n)-f_H(a))$ and $V=\frac{\sqrt{T_n}}{\sigma_H}E_2(n)$: by Theorem~\ref{thm:Q-cont} and Lemma~\ref{lem:quad},
$\dW\le C\nu_H(T_n)+C\sqrt{T_n}\,(T_n^{-1}+\Delta_n^{\gamma_H}T_n^{-1/2})\le C(\nu_H(T_n)+\Delta_n^{\gamma_H})$, because $T_n^{-1/2}\le\nu_H(T_n)$.

\emph{The estimator $\hat a_n$.} Write $\hat a_n-a=(f_H^{-1})'(\zeta_n)(\bQ_n-f_H(a))$ with $\zeta_n$ between $\bQ_n$ and $f_H(a)$, and $(f_H^{-1})'(\zeta_n)=g_H(a)+(f_H^{-1})''(\delta_n)(\zeta_n-f_H(a))$ with $\delta_n$ between $\bQ_n$ and $f_H(a)$. By the triangle inequality for $\dW$,
\[
\dW\Bigl(\frac{\sqrt{T_n}}{g_H\sigma_H}(\hat a_n-a),N\Bigr)\le\frac{\sqrt{T_n}}{|g_H|\sigma_H}\E\bigl[|(f_H^{-1})''(\delta_n)|\,|\bQ_n-f_H(a)|^2\bigr]+\dW\Bigl(\frac{\sqrt{T_n}}{\sigma_H}(\bQ_n-f_H(a)),N\Bigr).
\]
For the first term, $(f_H^{-1})''(x)=c\,x^{-r}$ with $r=2+\frac1{2H}$ is decreasing, so $|(f_H^{-1})''(\delta_n)|\le C+C\bQ_n^{-r}$, and by Hölder's inequality, Proposition~\ref{prop:neg-discrete} and hypercontractivity \eqref{eq:hypercontractivity}, we get :
\begin{align*}
\sqrt{T_n}\,\E\bigl[|(f_H^{-1})''(\delta_n)||\bQ_n-f_H|^2\bigr]&\le C\sqrt{T_n}\,\norm{\bQ_n-f_H(a)}{L^2}^2\\
&\le C\sqrt{T_n}\bigl(\tau_H(T_n)^2+\norm{E_2(n)}{2}^2\bigr)\\
&\le C\bigl(\nu_H(T_n)+T_n^{-3/2}+\Delta_n^{2\gamma_H}T_n^{-1/2}\bigr)\le C\nu_H(T_n).
\end{align*}
(ii) Apply \eqref{eq:coupling} with $U=T_n^{1-H}(\hat b_{T_n}-b)$, $V=T_n^{1-H}E_1(n)$: Theorem~\ref{thm:b-cont} gives the term $T_n^{-H}$ and Lemma~\ref{lem:linear} gives $T_n^{1-H}\norm{E_1}{1}\le CT^{1-H}\Delta_n^{\gamma_H}\tau_H(T_n)$.
\end{proof}

\begin{corollary}[Central limit theorems]\label{cor:CLT-discrete}
Let $\Delta_n\to0$ and $T_n\to\infty$.
\begin{enumerate}
\item[(i)] If $0<H<\frac34$, then $\sqrt{T_n}(\bQ_n-f_H(a))\Rightarrow\mathcal N(0,\sigma_H^2)$ and $\sqrt{T_n}(\hat a_n-a)\Rightarrow\mathcal N(0,g_H(a)^2\sigma_H^2)$.
\item[(ii)] If $H>\frac12$ then $T_n^{1-H}(\hat b_n-b)\Rightarrow\mathcal N(0,a^{-2})$ with no further condition; if $H\le\frac12$ the same holds provided that $T_n^{1-2H}\Delta_n^{1+2H}\to0$.
\end{enumerate}
\end{corollary}


\subsection{Numerical illustration}\label{subsec:numerics}

We illustrate Theorem~\ref{thm:BE-discrete} and Corollary~\ref{cor:CLT-discrete} by a short Monte Carlo study. We take $a=b=1$ and $X_0=0$. By the decomposition $X_t=b+Z_t^H-e^{-at}(b+Z_0^H)$, the vector $(Z^H_{t_i})_{0\le i<n}$ is a centered stationary Gaussian vector with covariance matrix $(\rho_H(|i-j|\Delta_n))_{i,j}$, so the discrete path of $X$ is simulated \emph{exactly}, as a Gaussian vector with this covariance matrix. The covariance is computed from the spectral representation
\[
\rho_H(u)=\int_{\mathbb R}e^{iu\lambda}\,\frac{\Gamma(2H+1)\sin(\pi H)}{2\pi}\,\frac{|\lambda|^{1-2H}}{a^2+\lambda^2}\,d\lambda ,
\]
by numerical quadrature. We check that : $\rho_H(0)=H\Gamma(2H)a^{-2H}$ and that $\rho_H(u)\sim H(2H-1)a^{-2}u^{2H-2}$ as $u\to\infty$ for $H\ne\frac12$.

\begin{figure}[ht]
\centering
\includegraphics[width=0.92\textwidth]{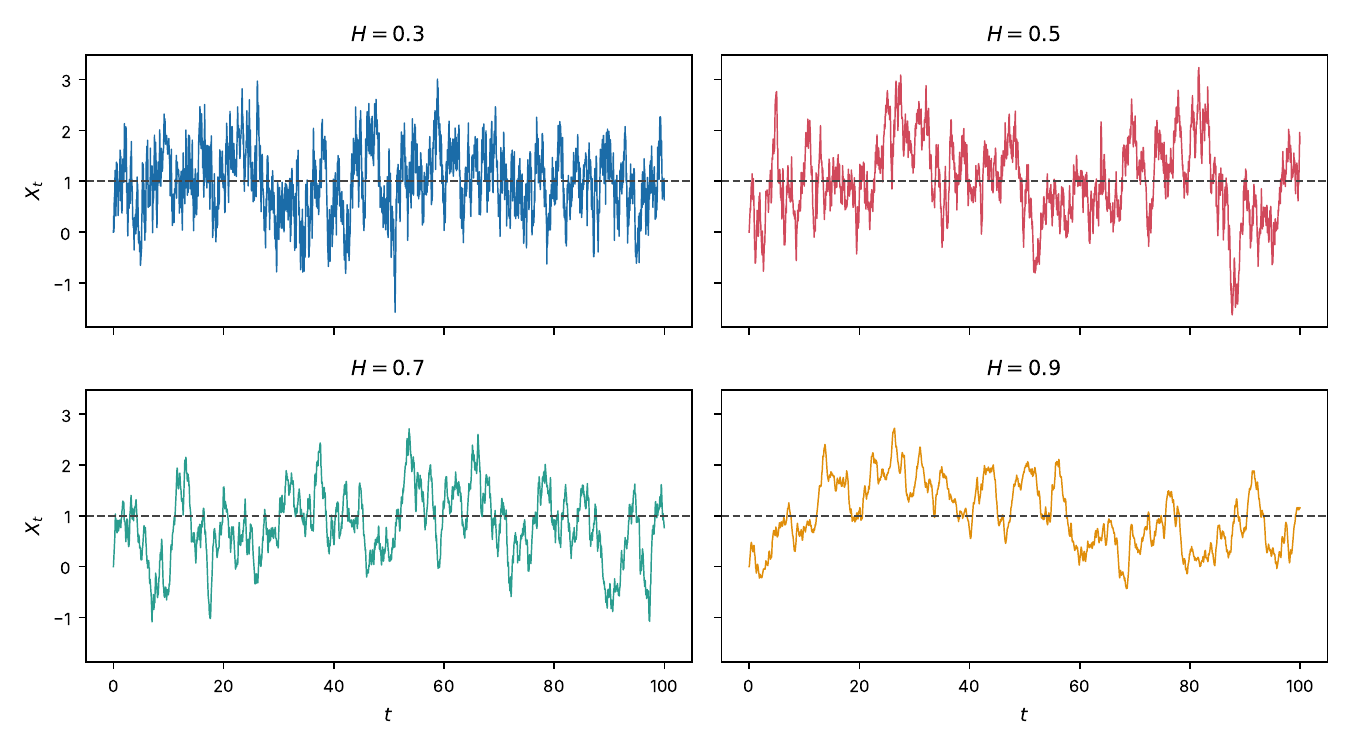}
\caption{One simulated trajectory of the fractional Vasicek process for $H=0.3,0.5,0.7,0.9$ ($a=b=1$, $X_0=0$, $n=2000$, $\Delta_n=0.05$, $T_n=100$). The dashed line is the level $b=1$.}\label{fig:traj}
\end{figure}

Figure~\ref{fig:traj} shows one simulated trajectory of $X$ for each of $H\in\{0.3,0.5,0.7,0.9\}$ ($a=b=1$, $X_0=0$, $n=2000$, $\Delta_n=0.05$, $T_n=100$). In every case the path starts at $0$, is quickly attracted to the level $b=1$ (dashed line) and then fluctuates around it, in accordance with the stationary behaviour of $X$. The regularity of the paths changes with the Hurst index: they are rough for $H=0.3$, as for a Brownian motion at $H=0.5$, and increasingly smooth with long excursions away from $b$ for $H=0.7$ and $0.9$, which is the long memory of the noise.

\begin{figure}[ht]
\centering
\includegraphics[width=0.92\textwidth]{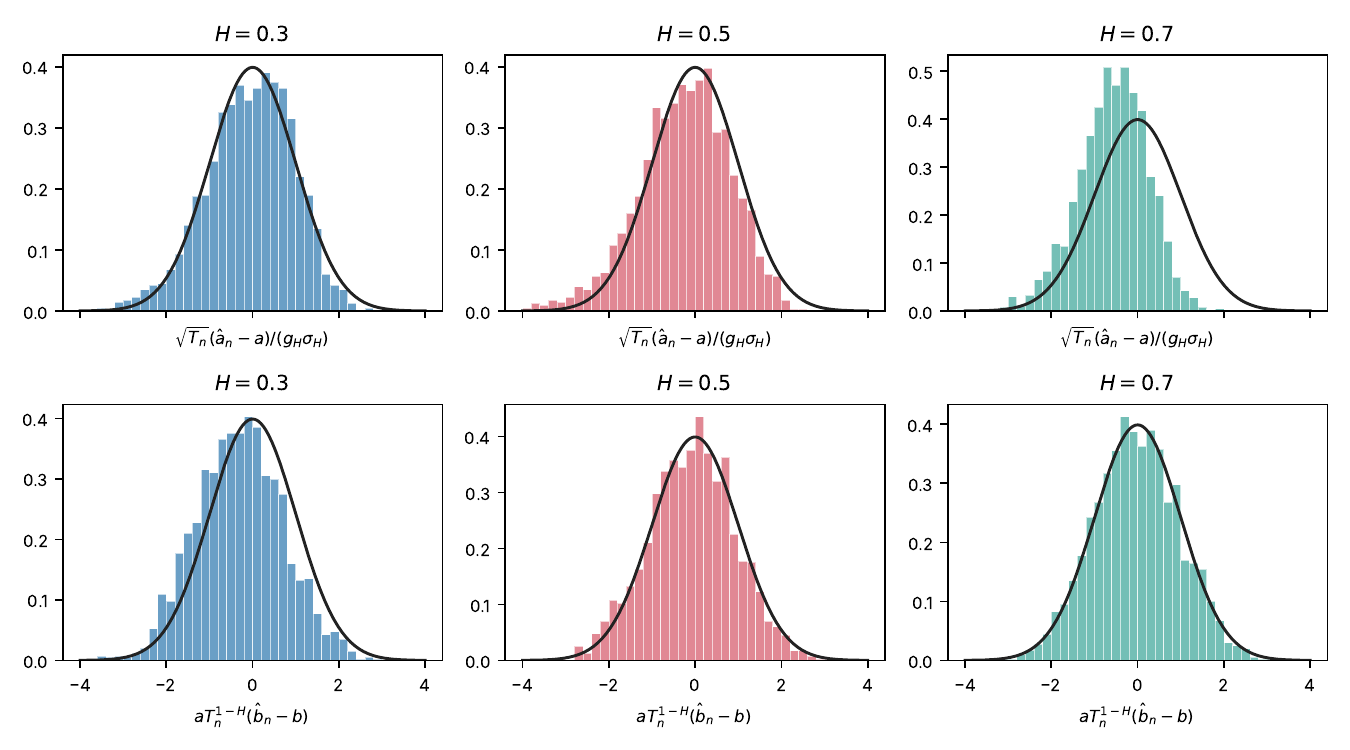}
\caption{Histograms of $\sqrt{T_n}(\hat a_n-a)/(g_H\sigma_H)$ (top) and $aT_n^{1-H}(\hat b_n-b)$ (bottom) for $H=0.3,0.5,0.7$ ($n=2000$, $\Delta_n=0.05$, $T_n=100$), with the $\mathcal N(0,1)$ density.}\label{fig:sim}
\end{figure}

We now turn to the statistics. The constant $\sigma_H^2=4\int_0^\infty\rho_H^2$ is obtained by numerical integration: $\sigma_H^2\approx0.169,\ 0.500,\ 2.941$ for $H=0.3,\,0.5,\,0.7$. We consider $H\in\{0.3,0.5,0.7\}$ and two sampling schemes with the same horizon $T_n=100$: $(n,\Delta_n)=(400,0.25)$ and $(2000,0.05)$. For each of $R=2000$ independent paths we compute $\hat a_n$ and $\hat b_n$ and the normalized statistics
\[
\frac{\sqrt{T_n}}{g_H(a)\,\sigma_H}\,(\hat a_n-a)\qquad\text{and}\qquad a\,T_n^{1-H}(\hat b_n-b),
\]
which are approximately $\mathcal N(0,1)$ by Corollary~\ref{cor:CLT-discrete}. Table~\ref{tab:sim} reports their empirical mean and standard deviation, and the empirical Wasserstein distance $d_W$ to $\mathcal N(0,1)$ (computed by the quantile coupling); Figure~\ref{fig:sim} shows the histograms for $\Delta_n=0.05$. The Monte Carlo error is about $0.02$ on the means and $0.02$ on the standard deviations.

\begin{table}[ht]
\centering
\small
\begin{tabular}{cc c | ccc | ccc}
\hline
 & & & \multicolumn{3}{c|}{$\hat a_n$} & \multicolumn{3}{c}{$\hat b_n$}\\
$H$ & $n$ & $\Delta_n$ & mean & s.d. & $d_W$ & mean & s.d. & $d_W$\\
\hline
0.3 & 400 & 0.25 & -0.07 & 1.10 & 0.09 & -0.25 & 1.04 & 0.25 \\
0.3 & 2000 & 0.05 & -0.06 & 1.01 & 0.07 & -0.30 & 1.00 & 0.30 \\
\hline
0.5 & 400 & 0.25 & -0.27 & 1.04 & 0.27 & -0.12 & 1.00 & 0.12 \\
0.5 & 2000 & 0.05 & -0.25 & 1.07 & 0.25 & -0.08 & 1.01 & 0.08 \\
\hline
0.7 & 400 & 0.25 & -0.55 & 0.80 & 0.55 & -0.06 & 0.99 & 0.06 \\
0.7 & 2000 & 0.05 & -0.57 & 0.81 & 0.57 & -0.03 & 0.99 & 0.04 \\
\hline
\end{tabular}
\caption{Normalized statistics for $\hat a_n$ and $\hat b_n$ ($a=b=1$, $T_n=100$, $R=2000$ replications): empirical mean, standard deviation and Wasserstein distance $d_W$ to $\mathcal N(0,1)$.}\label{tab:sim}
\end{table}

. For $\hat b_n$ the normalized statistic is close to $\mathcal N(0,1)$ for all $H$, with a standard deviation within $0.05$ of $1$; its negative mean is the bias due to the initial condition $X_0=0$: indeed $\E[\hat b_n-b]\approx-b/(aT_n)$, so that $\E[aT_n^{1-H}(\hat b_n-b)]\approx-b\,T_n^{-H}$, i.e.\ $-0.25,\,-0.10,\,-0.04$ for $H=0.3,\,0.5,\,0.7$, which agrees with the table up to the Monte Carlo error and is exactly the $T_n^{-H}$ term of the bound. For $\hat a_n$ the approximation is good for $H=0.3$ and $H=0.5$ (standard deviation close to $1$; the bias of order $T_n^{-1/2}$ is visible for $H=\frac12$), whereas for $H=0.7$, close to $\frac34$, the rate $\nu_H(T_n)=T_n^{2H-3/2}=T_n^{-0.1}$ is very slow and the standard deviation is still about $0.8$ at $T_n=100$. Finally, at fixed horizon, refining the grid from $\Delta_n=0.25$ to $0.05$ slightly improves the standard deviation for $H=0.3$ (from $1.10$ to $1.01$ for $\hat a_n$, from $1.04$ to $1.00$ for $\hat b_n$), in line with the discretization terms $\Delta_n^{H+1/2}$ and $T_n^{1/2-H}\Delta_n^{H+1/2}$, which are the largest for small $H$.

\section{Appendix}\label{sec:appendix}
\renewcommand{\thesubsection}{\Alph{subsection}}
\numberwithin{theorem}{subsection}

\subsection{Covariance of the stationary process}\label{app:cov}

Let $B^H$ be a fractional Brownian motion, $Z_t^H=\int_{-\infty}^te^{-a(t-s)}dB^H_s$, $\rho_H(u)=\E[Z_u^HZ_0^H]$ and $f_H(a)=\rho_H(0)$. Below, $C$ is a constant depending only on $(a,H)$.

\begin{lemma}\label{lem:rho}
\begin{enumerate}
\item[(a)] $|\rho_H(u)|\le C(1+u)^{2H-2}$ for $u\ge0$ (and $\rho_H(u)=e^{-au}/(2a)$ if $H=\frac12$); hence $\int_0^T|\rho_H(u)|du\le C\,T^{(2H-1)_+}$ for $T\ge1$;
\item[(b)] $0\le f_H-\rho_H(u)\le C u^{2H}$ for $u\ge0$;
\item[(c)] $\rho_H\in C^2(0,\infty)$; for $0<u\le1$, $|\rho_H'(u)|\le Cu^{2H-1}$ and $|\rho_H''(u)|\le Cu^{2H-2}$ ($\rho_H''$ is bounded if $H=\frac12$); for $u\ge1$, $|\rho_H'(u)|\le Cu^{2H-3}$ and $|\rho_H''(u)|\le Cu^{2H-4}$.
\end{enumerate}
Consequently, $|\rho_H(u)^2-\rho_H(0)^2|\le Cu^{2H}$ and $|(\rho_H^2)''(u)|\le C\,m(u)$, with $m(u):=u^{2H-2}\mathbf 1_{\{u\le1\}}+u^{4H-6}\mathbf 1_{\{u>1\}}$ (for $H=\frac12$, $(\rho_H^2)''$ is bounded on $(0,1]$).
\end{lemma}

\begin{proof}
If $H=\frac12$, It\^o's isometry gives $\rho_H(u)=e^{-au}/(2a)$ and everything is immediate. Let $H\ne\frac12$ and $\alpha_H=H(2H-1)$.

\emph{Representation of $\rho_H$.} Since $e^{as}$ is $C^1$ and $B^H$ is continuous, integration by parts holds pathwise (Riemann--Stieltjes) and the integral coincides a.s.\ with the Wiener integral. For $M>0$,
\[
\int_{-M}^0 e^{as}dB^H_s=e^{-aM}\bigl(B^H_0-B^H_{-M}\bigr)+a\int_0^M e^{-ar}\bigl(B^H_0-B^H_{-r}\bigr)dr .
\]
As $M\to\infty$ the boundary term tends to $0$ in $L^2$ (its norm is $e^{-aM}M^H$) and the integral converges in $L^2$ by Minkowski's inequality, so $Z_0^H=a\int_0^\infty e^{-ar}(B^H_0-B^H_{-r})\,dr$. Let $\mu(v):=\E[Z_0^H(B^H_v-B^H_0)]$, $v>0$. Since $\E[(B^H_0-B^H_{-r})(B^H_v-B^H_0)]=\frac12\bigl((v+r)^{2H}-r^{2H}-v^{2H}\bigr)$, Fubini's theorem gives $\mu(0+)=0$, and differentiating under the integral and integrating by parts in $r$,
\begin{gather*}
\mu'(v)=aH\!\int_0^\infty\! e^{-ar}\bigl[(v+r)^{2H-1}-v^{2H-1}\bigr]dr=\tilde q_H(v),\\
\tilde q_H(v):=\alpha_H\int_0^\infty e^{-az}(v+z)^{2H-2}dz .
\end{gather*}
Now write $Z_u^H=e^{-au}Z_0^H+\int_0^ue^{-a(u-v)}dB^H_v$, multiply by $Z_0^H$, take expectations and integrate by parts in $v$ (recall $\mu(0+)=0$):
\begin{equation}\label{eq:rho-rep}
\begin{gathered}
\rho_H(u)=e^{-au}f_H+\int_0^ue^{-aw}\tilde q_H(u-w)\,dw,\\
\rho_H'=-a\rho_H+\tilde q_H,\qquad \rho_H''=-a\rho_H'+\tilde q_H'\quad(u>0).
\end{gathered}
\end{equation}

\emph{Bounds on $\tilde q_H$.} The function $\tilde q_H$ is $C^\infty$ on $(0,\infty)$ and, by differentiation under the integral,
$|\tilde q_H(v)|\le Cv^{2H-2}$ for all $v>0$; $|\tilde q_H'(v)|\le Cv^{2H-2}$ for $0<v\le1$ (indeed $\tilde q_H'(v)=\alpha_H(2H-2)\int_0^\infty e^{-az}(v+z)^{2H-3}dz$ has modulus at most $|\alpha_H|\,v^{2H-2}$); $|\tilde q_H'(v)|\le Cv^{2H-3}$ for $v\ge1$; $|\tilde q_H''(v)|\le Cv^{2H-4}$ for $v\ge\frac12$. Near $0$: if $H<\frac12$, $|\tilde q_H(v)|\le Hv^{2H-1}$; if $H>\frac12$, $\tilde q_H(0+)=\alpha_H\Gamma(2H-1)a^{1-2H}=af_H$ and $|\tilde q_H(v)-\tilde q_H(0+)|\le Hv^{2H-1}$. In all cases $\int_0^x|\tilde q_H|\le C(1+x)$.
\begin{enumerate}
	\item \emph{(a)}  For $u\ge1$ split $\int_0^ue^{-aw}\tilde q_H(u-w)dw$ at $w=u/2$: on $w\le u/2$, $|\tilde q_H(u-w)|\le Cu^{2H-2}$; on $w>u/2$ the integral is at most $e^{-au/2}\int_0^{u}|\tilde q_H|\le Ce^{-au/2}(1+u)$. Hence $|\rho_H(u)|\le Cu^{2H-2}$, and $|\rho_H|\le f_H$ for $u\le1$. The bound on $\int_0^T|\rho_H|$ follows by integrating $(1+u)^{2H-2}$.
	\item \emph{(b)} $2(f_H-\rho_H(u))=\E(Z_u^H-Z_0^H)^2\ge0$ by stationarity, and $Z_u^H-Z_0^H=-a\int_0^uZ_v^Hdv+B^H_u-B^H_0$ gives $\|Z_u^H-Z_0^H\|_{L^2}\le au\sqrt{f_H}+u^H\le Cu^H$ for $u\le1$; for $u\ge1$, $f_H-\rho_H\le2f_H\le2f_Hu^{2H}$.
	
	\item \emph{(c)} Let $0<u\le1$. If $H>\frac12$, then $\rho_H'(u)=-a(\rho_H(u)-f_H)+(\tilde q_H(u)-\tilde q_H(0+))$, so by (b) $|\rho_H'(u)|\le Cu^{2H}+Hu^{2H-1}\le Cu^{2H-1}$. If $H<\frac12$, then $|\rho_H'|\le a|\rho_H|+|\tilde q_H|\le af_H+Hu^{2H-1}\le Cu^{2H-1}$. In both cases $|\rho_H''|\le a|\rho_H'|+|\tilde q_H'|\le Cu^{2H-1}+Cu^{2H-2}\le Cu^{2H-2}$.
	Let $u\ge1$. Since $\tilde q_H(u)=a\int_0^ue^{-aw}\tilde q_H(u)\,dw+e^{-au}\tilde q_H(u)$, \eqref{eq:rho-rep} gives
	\[
	\rho_H'(u)=e^{-au}\bigl(\tilde q_H(u)-af_H\bigr)+a\int_0^ue^{-aw}\bigl[\tilde q_H(u)-\tilde q_H(u-w)\bigr]dw .
	\]
	For $w\le u/2$ use $\tilde q_H(u)-\tilde q_H(u-w)=w\tilde q_H'(u)+O(w^2u^{2H-4})$; the range $w>u/2$ and the first term are $O(e^{-au/2}(1+u))$. Therefore $\rho_H'(u)=\tilde q_H'(u)/a+\mathcal E(u)$ with $|\mathcal E(u)|\le Cu^{2H-4}$, so $|\rho_H'|\le Cu^{2H-3}$ and $\rho_H''=-a\rho_H'+\tilde q_H'=-a\,\mathcal E$ satisfies $|\rho_H''|\le Cu^{2H-4}$.
	\item \emph{Consequences for $\rho_H^2$.} $\rho_H(0)^2-\rho_H(u)^2=(f_H-\rho_H)(f_H+\rho_H)$ gives the first bound by (b). Finally $(\rho_H^2)''=2\rho_H'^2+2\rho_H\rho_H''$: for $u\le1$, $\rho_H'^2\le Cu^{4H-2}\le Cu^{2H-2}$ and $|\rho_H\rho_H''|\le Cu^{2H-2}$; for $u\ge1$, $\rho_H'^2$ and $|\rho_H\rho_H''|$ are $\le Cu^{4H-6}$.
	
\end{enumerate}
\end{proof}

\subsection{Negative moments of $\bQ(T)$ and $\bQ_n$}\label{app:neg}

We will make use of the Mandelbrot--Van Ness representation \cite{MandelbrotVanNess1968}, recalled in \eqref{eq:MVN}, of the fractional Brownian motion: Let $W$ be a two-sided standard Brownian motion and
\[
B^H_t=\int_{\R}K(t,u)\,dW_u,\qquad K(t,u)=c_H\bigl[(t-u)_+^{H-1/2}-(-u)_+^{H-1/2}\bigr],
\]
where $c_H>0$ normalizes $\Var B^H_1=1$. Let $\mathcal F_t=\sigma(W_u-W_v:\,v\le u\le t)$. Then $B^H_s$, $s\le t$, is $\mathcal F_t$-measurable, and $W$-integrals over $(t,\infty)$ are independent of $\mathcal F_t$. We write $X_t=k(t)+G_t$, $k(t)=b(1-e^{-at})$, $G_t=\int_0^te^{-a(t-s)}dB^H_s$.

\subsubsection{The innovation kernel}

Integrating by parts (in the Riemann–Stieltjes sense), $G_t=B^H_t-a\int_0^te^{-a(t-s)}B^H_sds$, and stochastic Fubini's theorem gives
\begin{equation}\label{eq:Gkernel}
G_t=\int_{\R}\Bigl[K(t,u)-a\!\int_0^te^{-a(t-s)}K(s,u)\,ds\Bigr]dW_u .
\end{equation}
Fix $a_0\ge0$ and $t\ge a_0$. For $u>a_0\ge0$ we have $(-u)_+=0$ and $K(s,u)=c_H(s-u)_+^{H-1/2}$, so the part of \eqref{eq:Gkernel} coming from $u\in(a_0,t]$ is
\begin{equation}\label{eq:innovation}
A_t:=\int_{a_0}^t\kappa(t-u)\,dW_u,\qquad
\kappa(\tau):=c_H\Bigl[\tau^{H-1/2}-a\int_0^\tau e^{-a(\tau-r)}r^{H-1/2}dr\Bigr],
\end{equation}
which is independent of $\mathcal F_{a_0}$, while $C_t:=G_t-A_t$ is $\mathcal F_{a_0}$-measurable.

\begin{lemma}\label{lem:kappa}
Let $\theta_0:=\frac{H+1/2}{2a}$ and $q(v):=\int_0^v\kappa(\tau)^2d\tau$. Then $q(v)\ge c_1\,\min\{v,\theta_0\}^{2H}$ for all $v\ge0$, with $c_1=c_H^2/(8H)$. Moreover, for $a_0\le s\le t$,
\[
\E(A_t-A_s)^2\ge q(t-s).
\]
\end{lemma}

\begin{proof}
Since $0\le a\int_0^\tau e^{-a(\tau-r)}r^{H-1/2}dr\le\frac{a\tau^{H+1/2}}{H+1/2}\le\frac12\tau^{H-1/2}$ for $\tau\le\theta_0$, we get $\kappa(\tau)\ge\frac{c_H}2\tau^{H-1/2}$ there, hence $q(v)\ge\frac{c_H^2}{4}\int_0^{v\wedge\theta_0}\tau^{2H-1}d\tau=c_1(v\wedge\theta_0)^{2H}$ ($q$ is non-decreasing). For the second claim, $A_t-A_s=\int_{a_0}^s(\kappa(t-u)-\kappa(s-u))dW_u+\int_s^t\kappa(t-u)dW_u$, the two integrals are independent and the variance of the second is $q(t-s)$.
\end{proof}

For $H=\frac12$ one has $\kappa(\tau)=e^{-a\tau}$ and $q(v)=(1-e^{-2av})/(2a)$.

\subsubsection{Continuous time}

\begin{proposition}\label{prop:neg-cont}
For every $p\ge1$ there are $T_0$ and $C_p$ such that $\E[\bQ(T)^{-p}]\le C_p$ for all $T\ge T_0$.
\end{proposition}

\begin{proof}
Fix an integer $m>4p$ and cut $[0,T]$ into $m$ blocks $I_k=[a_k,a_k+L)$, $a_k=(k-1)L$, $L=T/m$. Set $\bar X_{I_k}=\frac1L\int_{I_k}X_t\,dt$ and $\sigma_k^2=\frac1L\int_{I_k}(X_t-\bar X_{I_k})^2dt$.
\\
\\
For each block, $\int_{I_k}(X_t-c)^2dt=\int_{I_k}(X_t-\bar X_{I_k})^2dt+L(\bar X_{I_k}-c)^2$; taking $c=\bar X_T$ and using the AM--GM inequality,
\[
\bQ(T)\ge\frac1m\sum_{k=1}^m\sigma_k^2\ge\Bigl(\prod_{k=1}^m\sigma_k^2\Bigr)^{1/m},\qquad\text{so}\qquad \bQ(T)^{-p}\le\prod_{k=1}^m(\sigma_k^2)^{-p/m}.
\]
\\
\\
By \eqref{eq:innovation} with $a_0=a_k$, for $t\in I_k$, $X_t=A^k_t+D^k_t$ where $A^k$ is independent of $\mathcal F_{a_k}$ and $D^k_t=k(t)+C^k_t$ is $\mathcal F_{a_k}$-measurable. We have,
\[
\E[\sigma_k^2\mid\mathcal F_{a_k}]\ \ge\ \frac1L\int_{I_k}\E(A^k_t-\bar A^k_{I_k})^2dt=\frac1{2L^2}\iint_{I_k^2}\E(A^k_t-A^k_s)^2\,ds\,dt,
\]
where we used the identity $\frac1L\int(U_t-\bar U)^2dt=\frac1{2L^2}\iint(U_t-U_s)^2dsdt$ .
\\
\\
By Lemma~\ref{lem:kappa}, and since $\{|t-s|\ge\theta_0\}$ has Lebesgue measure $(L-\theta_0)^2\ge L^2/4$ in $I_k^2$ when $L\ge2\theta_0$,
\[
\E[\sigma_k^2\mid\mathcal F_{a_k}]\ge\frac{c_1\theta_0^{2H}}{8}=:c_2\qquad\text{whenever }T\ge T_0:=2m\theta_0 .
\]
Note that $T_0$ depends only on $m$, $a$ and $H$.
\\
\\
Applying the Carbery--Wright inequality \cite{CarberyWright2001}, recalled in \eqref{eq:CW}, there is a universal $c>0$ with
\[
\Prob\bigl(\sigma_k^2\le\varepsilon\mid\mathcal F_{a_k}\bigr)\le c\sqrt{\varepsilon/c_2}\qquad(\varepsilon>0,\ T\ge T_0).
\]
Therefore
\begin{align*}
\E\bigl[(\sigma_k^2)^{-p/m}\mid\mathcal F_{a_k}\bigr]&=\int_0^\infty\Prob\bigl(\sigma_k^2<x^{-m/p}\mid\mathcal F_{a_k}\bigr)dx\\
&\le1+\frac c{\sqrt{c_2}}\int_1^\infty x^{-m/(2p)}dx=1+\frac c{\sqrt{c_2}}\,\frac1{m/(2p)-1}=:\gamma_m<\infty,
\end{align*}
the integral converging because $m>4p>2p$.
\\
\\
For $j<k$, $\sigma_j^2$ depends on $X_t$, $t\le a_k$, hence is $\mathcal F_{a_k}$-measurable. Conditioning successively for $k=m,m-1,\dots,1$,
$\E[\bQ(T)^{-p}]\le\E\prod_k(\sigma_k^2)^{-p/m}\le\gamma_m^{\,m}$ for all $T\ge T_0$.
\end{proof}

\subsubsection{Discrete time (proof of Proposition~\ref{prop:neg-discrete})}

The proof follows that of Proposition~\ref{prop:neg-cont}, replacing integrals by sums over the grid; we detail the steps that differ.

Fix $p\ge1$ and an integer $m>4p$, and let $n\ge2m$, $\ell:=\lfloor n/m\rfloor$. Define the blocks of $\ell$ consecutive grid points $\mathcal I_k=\{t_i:\,(k-1)\ell\le i<k\ell\}$, $k=1,\dots,m$, with $a_k=t_{(k-1)\ell}$, and
$\bar X_{\mathcal I_k}=\frac1\ell\sum_{\mathcal I_k}X_{t_i}$, $s_k^2=\frac1\ell\sum_{\mathcal I_k}(X_{t_i}-\bar X_{\mathcal I_k})^2$.
\\
\\
Dropping the (non-negative) terms with $i\ge m\ell$ and using $\sum_{\mathcal I_k}(X_{t_i}-c)^2\ge\sum_{\mathcal I_k}(X_{t_i}-\bar X_{\mathcal I_k})^2$ with $c=\hat b_n$,
\[
\bQ_n\ge\frac{\ell}{n}\sum_{k=1}^ms_k^2\ge\frac1{2m}\sum_{k=1}^ms_k^2\ge\frac12\Bigl(\prod_ks_k^2\Bigr)^{1/m},
\]
since $\ell\ge n/m-1\ge n/(2m)$. Thus $\bQ_n^{-p}\le2^p\prod_k(s_k^2)^{-p/m}$.
\\
\\
With $A^k$, $D^k$ as in \eqref{eq:innovation} and the discrete variance identity $\frac1\ell\sum_i(U_i-\bar U)^2=\frac1{2\ell^2}\sum_{i,j}(U_i-U_j)^2$, Lemma~\ref{lem:kappa} gives
\[
\E[s_k^2\mid\mathcal F_{a_k}]\ \ge\ \frac1{2\ell^2}\sum_{t_i,t_j\in\mathcal I_k}q(|t_i-t_j|)\ \ge\ \frac{c_1\theta_0^{2H}}{2\ell^2}\,\#\{(i,j):\,|i-j|\Delta_n\ge\theta_0\}.
\]
The number of pairs with $|i-j|\Delta_n<\theta_0$ is at most $\ell(2\theta_0/\Delta_n+1)$, so the right-hand side is $\ge\frac{c_1\theta_0^{2H}}{4}=:c_3$ as soon as $\ell\ge2(2\theta_0/\Delta_n+1)$, i.e.\ as soon as the block length $\ell\Delta_n$ satisfies $\ell\Delta_n\ge4\theta_0+2\Delta_n$. This holds for $n\ge n_0(a,H,m)$ because $T_n/m\to\infty$ and $\Delta_n\le1$ (precisely, it suffices that $T_n\ge m(4\theta_0+4)$).
\\
\\
Conditionally on $\mathcal F_{a_k}$, $s_k^2$ is a degree-$2$ polynomial of a Gaussian vector with conditional mean $\ge c_3$. The Carbery--Wright inequality \cite{CarberyWright2001} gives $\Prob(s_k^2\le\varepsilon\mid\mathcal F_{a_k})\le c\sqrt{\varepsilon/c_3}$, hence $\E[(s_k^2)^{-p/m}\mid\mathcal F_{a_k}]\le\gamma_m':=1+\frac{c}{\sqrt{c_3}}\frac1{m/(2p)-1}$, uniformly in $n\ge n_0$ and $k$. As $\prod_{j<k}(s_j^2)^{-p/m}$ is $\mathcal F_{a_k}$-measurable, iterating the conditioning for $k=m,\dots,1$ yields $\E[\bQ_n^{-p}]\le2^p(\gamma_m')^m$ for all $n\ge n_0$. \qed

\end{document}